\documentclass[a4,10pt]{article}
\usepackage{graphicx}
\usepackage{amsmath,amssymb}
\usepackage{ulem}

\usepackage{color}
\usepackage{here}

\begin{document}
\setlength{\baselineskip}{15pt} 
\def\pot{\mathaccent"7017} 
\makeatletter
 \renewcommand{\theequation}{%
  \thesection.\arabic{equation}}
 \@addtoreset{equation}{section}
\makeatother

\newtheorem{remark}{Remark}[section]
\newtheorem{theo}{Theorem}[section]
\newtheorem{lemma}{Lemma}[section]
\newtheorem{prop}{Proposition}[section]
\newtheorem{assume}{Assumption}[section]
\newtheorem{cor}{Corollary}[section]
\newtheorem{example}{Example}[section]

\def\ep{\varepsilon}
\def\ov{\overline}
\def\un{\underline}
\def\del{\partial}
\def\norm{\parallel}
\def\no{\noindent}
\def\lam{\lambda}
\def\dis{\displaystyle}
\def\bhat{\widehat}
\def\lap{\bigtriangleup}
\def\R{\mbox{\bf R}}
\def\C{\mbox{\bf C}}
\def\lg{\langle}
\def\rg{\rangle}

\newcommand{\Qed}{ $\square$ }

\title{
Single Transition Layer in Mass-Conserving Reaction-Diffusion 
Systems with Bistable Nonlinearity
}

\author{
Masataka Kuwamura\thanks{
Corresponding author; 
Graduate School of Human Development and Environment,
Kobe University, Kobe 
657-8501, Japan,
email: kuwamura@main.h.kobe-u.ac.jp }\;
,
Takashi Teramoto\thanks{
Department of Data Science, Kyoto Women's University, Kyoto 605-8501, Japan, 
email: teramotot@kyoto-wu.ac.jp}\;
and
Hideo Ikeda\thanks{
Department of Mathematics, University of Toyama, Toyama,
930-8555, Japan,
email: hideoikeda5@gmail.com}
}

\date{}


\maketitle
\begin{abstract}
Mass-conserving reaction-diffusion 
systems with bistable nonlinearity
are useful models for studying cell polarity formation,
which is a key process in cell division and differentiation.
We rigorously show the existence and stability of stationary solutions with a
single internal transition layer in such reaction-diffusion systems
under general assumptions by the singular perturbation theory.
Moreover, we present a meaningful model for understanding the 
existence of an unstable transition layer solution; our numerical simulations show that 
the unstable solution is a separatrix of the dynamics of the model.
\end{abstract}

\vspace{1ex}

Abbreviated title:  Transition Layer in Reaction-Diffusion Systems

\vspace{1ex}

\hspace{2.3cm} transition layer

\vspace{1ex}


\section{Introduction } \label{intoro}

Reaction-diffusion systems are useful models for studying the mechanism of the appearance of non-uniform patterns in various fields of science and technology.
In this paper, we consider the existence and stability of stationary solutions with a
single internal transition layer in the following reaction-diffusion system:
\begin{equation}\label{a1}
\left \{
\begin{array}{rcl}
u_t & = & \ep^2 u_{xx} +  f(u,v)  \\[1ex]
v_t & = & Dv_{xx}  -  f(u,v)
\end{array} 
\right.
\end{equation}
on an interval $0 < x < 1$ under the Neumann boundary condition, 
where $\ep$ and $D$ are positive constants satisfying $0 < \ep \ll D$.
The nonlinear term $f$ is a smooth cubic function such that 
the ODE $u_t = f(u, v)$ is bistable in $u$ for each fixed $v$.
This type of reaction-diffusion system is proposed by \cite{MJE1} 
for studying the wave-pinning phenomenon in cell division and differentiation:
A transient and localized stimulus to 
an unpolarized cell is spatially amplified to result in a robust subdivision of the cell into
two clearly defined regions, front and back. 
In simple terms, this biological phenomenon can be mathematically interpreted as 
the dynamics of \eqref{a1} where a traveling front solution starting from the edge 
converges to a stationary solution with a single internal transition layer. 
In fact, \cite{MJE1,MJE2} concluded that 
\begin{equation}\label{a1x}
\left \{
\begin{array}{rcl}
u_{\tau} & = & \ep^2 u_{xx} +  f(u,v)  \\[1ex]
v_{\tau} & = & Dv_{xx}  -  f(u,v)
\end{array} 
\right.
\end{equation}
has stable stationary solutions with a single internal transition layer 
under certain conditions by using a formal analysis and a perturbative argument
against the background of cell biology.
Here, $\tau = t/\ep$ represents the fast time scale as compared to the time scale $t$, 
and the dynamics of \eqref{a1x} is equivalent to that of \eqref{a1}.
They confirmed their theoretical results by numerical simulations for \eqref{a1x}
with specific nonlinear terms.

We note that \eqref{a1} is a typical example of reaction-diffusion systems
\begin{equation}\label{a2xx}
\left \{
\begin{array}{rcl}
u_t & = & d_1 u_{xx} +  f(u,v)  \\[1ex]
v_t & = &d_2 v_{xx}  -  f(u,v),
\end{array} 
\right.
\end{equation}
where the nonlinear term $f$ is an appropriate smooth function. 
They were proposed by \cite{Is, Ot} as a conceptual model for studying cell polarity formation which plays a key role in cell division and differentiation.
We immediately find that
any (smooth) solution of \eqref{a2xx} satisfies
\[
\int_0^1 \left( u(x, 0) + v(x, 0) \right) dx \equiv 
\int_0^1 \left( u(x, t) + v(x, t) \right) dx
\]
under the Neumann or periodic boundary conditions.
Therefore, \eqref{a2xx} are called the mass-conserving reaction-diffusion systems
\cite{BJF}.
Similarly, we immediately find that 
any (smooth) stationary solution of \eqref{a2xx} 
satisfies 
\begin{equation}\label{a2zzz}
d_1 u(x) + d_2 v(x) \equiv Const.
\end{equation}
under the Neumann or periodic boundary conditions. This fact suggests that
the asymptotic behavior of \eqref{a2xx} near stationary solutions can be similar to that of
scalar reaction-diffusion equations.

We should consider that  \eqref{a1} or \eqref{a1x} can include the dynamics which is 
similar to the dynamics of scalar 
reaction-diffusion equations with cubic nonlinearity such as the Allen-Cahn equation.
Recalling that the Allen-Cahn equation exhibits metastable patterns with
a single internal transition layer moving at the speed of $O(e^{-C/\ep})$ for some positive constant $C$
as shown in \cite{CP, HF, KO},
we wonder whether or not numerical solutions regarded as
stable stationary solutions in \cite{MJE1, MJE2} are metastable patterns.
In fact, metastable patterns are often misidentified as stable stationary solutions 
in numerical simulations; this would be a problem
in the study of the bifurcation structure of 
mass conserving reaction-diffusion systems with bistable nonlinearity
referring to numerical simulations as seen in \cite{MJE2}.
Moreover, the formal analysis and perturbative argument by \cite{MJE1,MJE2}
cannot distinguish stable stationary solutions from
metastable patterns.
Therefore, it is necessary and important to give a rigorous proof of the assertion that
\eqref{a1} or \eqref{a1x} has stable stationary solutions with a single internal transition layer.

The purpose of this paper is to show the existence and stability 
of stationary solutions with a
single internal transition layer in reaction-diffusion 
systems \eqref{a1} under general assumptions by the singular perturbation method \cite{HS, I, MTH, NF}.
Our approach provides a new point of view for studying mass-conserving 
reaction-diffusion systems.
In fact, to the best of 
our knowledge, existing mathematical analyses for mass-conserving 
reaction-diffusion systems
are concerned with localized unimodal patterns (spike solutions) under the 
influence of the ideas of the variational method \cite{CMS,ET,JM,EKS,LS,MO}.
First, we mention assumptions to obtain our main result.

\begin{assume}\label{ass1}
{\rm (A1)} 
The ODE $u_t = f(u, v)$
is bistable in $u$ for each fixed $v \in I = (\un{v}, \ov{v})$.
That is, $f(u, v)=0$ has exactly three roots 
$h^-(v) < h^0(v) < h^+(v) $ for each $v \in I$ satisfying
\[
f_u ( h^{\pm}(v), v ) < 0
\ \ \
\text{and}
\ \ \ 
f_u ( h^{0}(v), v ) > 0. 
\]
{\rm (A2)} The function 
\begin{equation}\label{a4}
J(v) = \int_{h^-(v)}^{h^+(v)} f(u, v)du \ \ \ ( v \in I)
\end{equation}
has an isolated zero at $v = v^* \in I$ such that 
\begin{equation}\label{a4x}
J'(v^*) = \int_{h^-(v^*)}^{h^+(v^*)} f_v(u, v^*)du \neq 0.
\end{equation}

\no
{\rm (A3)} 
\[
f_u ( h^{\pm}(v), v ) < f_v ( h^{\pm}(v), v )  \ \ \ ( v \in I).
\]
\end{assume}

The assumptions (A1), (A2) and (A3) were used in the formal analysis and
perturbative argument by \cite{MJE2}. 
We emphasize that these assumptions can be easily verified in 
many practical problems including those of \cite{MJE1, MJE2}. 
Let
\[
\xi := \int_0^1 \left( u(x, 0) + v(x, 0) \right) dx.
\]
Then, the conservation law
\begin{equation}\label{a2}
\int_0^1 \left( u(x, t) + v(x, t) \right) dx \equiv \xi
\end{equation}
holds. We choose $\xi$ in such a way that the following inequality holds:
\begin{equation}\label{a6}
h^-(v^*) + v^* < \xi < h^+(v^*) + v^*.
\end{equation}

Let us define 
\begin{equation}\label{a8}
U^*(x) = 
\left \{
\begin{array}{l}
h^-(v^*)  \ \ \ (0 \leq x \leq x^*) \\[1ex]
h^+(v^*)  \ \ \ (x^* < x \leq 1) 
\end{array} 
\right.
\end{equation}
and
\begin{equation}\label{a9}
V^*(x) = v^* \ \ \ ( 0 \leq x \leq 1),
\end{equation}
where 
\begin{equation}\label{a7}
x^* := \dis\frac{ h^+(v^*) + v^* - \xi}{h^+(v^*) - h^-(v^*)}
\end{equation}
satisfies $0 < x^* < 1$ by \eqref{a6}. We note that \eqref{a7} is derived from
$$
\int_0^1 \left( U^*(x) + V^*(x) \right) dx = \xi,
$$
and that $(U^*(x), V^*(x))$ is 
the leading term of the outer approximation of stationary solutions 
constructed in Section~2.

\begin{figure}[ht!]
\centering
\includegraphics[width=12cm, bb=0 0 810 296]{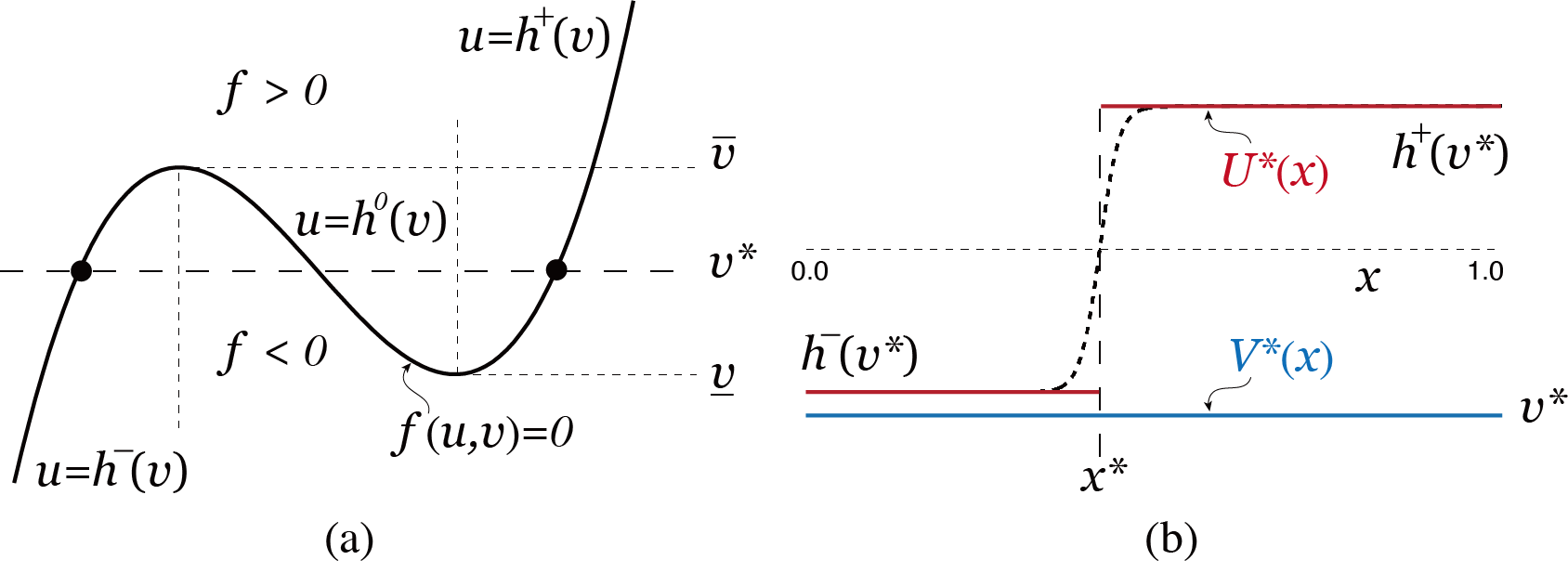}
  \caption{Schematic pictures of bistable nonlinearity and stationary solutions with a single internal transition layer
  at the leading order. (a) The graph of $f(u,v) = 0$ separates the $(u,v)$-plane into the upper-left and lower-right regions of $f > 0$ and $f<0$, respectively. 
Solid points on the graph indicate bistable equilibria of the ODE $u_t = f(u, v)$ in $u$
for $v = v^*$. (b) The red and blue solid lines show the profiles of $U^*(x)$ and $V^*(x)$, respectively.  
The dotted line shows the profile of $u^\ep(x)$ around $x = x^*$, 
which jumps up from $h^-(v^*)$ to $h^+(v^*)$ at $x = x^*$
given by \eqref{a7}.
  }
  \label{example0}
\end{figure}
In simple terms, our main result is summarized as follows:

\begin{theo}\label{th1}
Assume (A1) and (A2).
For any given $\xi$ satisfying \eqref{a6}, 
the reaction-diffusion system \eqref{a1} has a 
family of stationary solutions
$(u^{\ep}(x), v^{\ep}(x))$
satisfying 
\begin{equation}\label{a2xyz}
\int_0^1 \left( u^\ep(x) + v^\ep(x) \right) dx = \xi 
\end{equation}
for small $\ep$. They satisfy 
$$
\lim_{\ep \to 0} u^\ep(x)  =  U^*(x) \ \mbox{uniformly on} \ [0, x^* - \sigma]
\cup [x^* + \sigma, 1]
$$
for any $\sigma$ with $0 < \sigma < \min(x^*, 1-x^*)$, and 
$$
 \lim_{\ep \to 0} v^\ep(x)  =  V^*(x) \ \, \mbox{uniformly on} \ [0,1]. 
$$
Moreover, under the assumptions (A1), (A2) and (A3), the stationary solutions
$(u^{\ep}(x), v^{\ep}(x))$ are stable if $J'(v^*) > 0 $.
\end{theo}

\begin{remark}\label{rem0x} \rm
As shown in Figure~\ref{example0}(b), $u^\ep(x)$
exhibits
a single internal transition layer with $O(1)$-amplitude at $x = x^*$ when $\ep$ is sufficiently small. We note that it exhibits a jump-up transition. If we take 
$$
U^*(x) = 
\left \{
\begin{array}{l}
h^+(v^*)  \ \ \ (0 \leq x \leq x^*) \\[1ex]
h^-(v^*)  \ \ \ (x^* < x \leq 1),
\end{array} 
\right.
$$
we can obtain a stationary solution whose $u$-component exhibits a 
jump-down transition when $\ep$ is sufficiently small. 
It is stable if $J'(v^*) > 0 $.
\end{remark}

\begin{remark}\label{rem0xy} \rm
Noting \eqref{a2zzz} with \eqref{a2xx}, 
we immediately find $\ep^2 u^\ep(x) + D v^\ep(x) \equiv C(\ep)$,
where $C(\ep)$ denotes a constant independent of $x$. Therefore, 
we have $v^\ep(x) = (C(\ep) - \ep^2 u^\ep(x))/D$, which 
implies that $v^\ep(x)$ exhibits
a single internal transition layer with $O(\ep^2)$-amplitude at $x = x^*$
when $\ep$ is sufficiently small (See, Figures~\ref{example1}-\ref{example3} in 
Section~4).
Moreover, it
exhibits a jump-down (resp. jump-up) transition if $u^\ep(x)$ exhibits 
a jump-up (resp. jump-down) transition.
\end{remark}

The mathematically precise version of this theorem is presented by Theorem~\ref{th2},
Corollary~\ref{cor1} and Theorem~\ref{th4}.
They give a mathematical justification for the results in 
\cite{MJE1, MJE2}, which are obtained by a formal analysis and 
a perturbative argument with the aid of numerical simulations.

The organization of this paper is as follows.
In the next section, we construct stationary solutions with a single internal transition layer
by the matched asymptotic expansion method \cite{I} which was used for
scalar reaction-diffusion equations with bistable nonlinearity dependent on
the spatial variable $x$. 
The first key of our construction is to find a rule for recursively performing
the matched asymptotic expansion of any order under the constraint condition \eqref{a2xyz}.
Moreover, the second key is to apply the implicit function
theorem to guarantee that the stationary solutions exactly satisfy the constraint condition \eqref{a2xyz}. 
In Section \ref{stability}, we show the stability of the stationary solutions with a single internal
transition layer on the basis of the singular perturbation method by \cite{HS, NF} which were used for
reaction-diffusion systems with bistable nonlinearity of FitzHugh-Nagumo type.
Here, we cannot apply 
the Lax-Milgram theorem
to the singular limit eigenvalue problem (SLEP) for investigating the behavior of
a critical eigenvalue which essentially determines
the stability of the stationary solutions. 
In spite of this complication, we can give the precise characterization of the critical eigenvalue due to a natural constrained condition 
derived from the conservation law.
Consequently, we can show the stability of the 
stationary solutions with a single internal transition layer under 
Assumption~\ref{ass1} when $J'(v^*) > 0$.
In Section 4, we present some examples of mass-conserving reaction-diffusion systems 
which have
a stationary solution with a single internal transition layer 
with the aid of numerical simulations.
We show a precise profile of the stationary solution, where
the $v$-component also exhibits a transition layer
at $x = x^*$ as well as the $u$-component.
Moreover, our numerical simulations show that the stationary solution is unstable if $J'(v^*) < 0$, which supports that 
the sign of $J'(v^*)$ gives the stability criterion of
stationary solutions with a single internal transition layer.
Section \ref{conclude} is devoted to concluding remarks.

\section{Existence of single transition layer solutions } \label{existence} 

\noindent
In this section, we consider single transition layer solutions of 
\begin{equation}\label{b1}
\left \{
\begin{array}{l}
\begin{array}{l}
\ep^2 u_{xx} +  f(u,v) = 0, \\[1ex]
Dv_{xx}  - f(u, v) = 0, 
\end{array} \quad  x \in (0,1) \\[0.5cm]
(u_x, v_x)(0) = (u_x, v_x)(1) = (0, 0)
\end{array} 
\right.
\end{equation}
satisfying 
\begin{equation}\label{b1_1}
\xi = \int_0^1 \left( u(x) + v(x) \right) dx 
\end{equation}
for a given constant $\xi$ with \eqref{a6} 
under the assumptions (A1) and (A2).

Since $\ep^2 u_{xx} + Dv_{xx} = 0$ by \eqref{b1}, noting the Neumann
boundary condition, we have
\begin{equation}\label{b2zz}
\ep^2 u + Dv = C(\ep),
\end{equation}
where $C(\ep)$ is a constant dependent on $\ep$. 
Substituting $v = ( C(\ep) - \ep^2 u )/D$ into the first equation of \eqref{b1}, 
and \eqref{b1_1}, we have a single equation for $u$ 
\begin{equation}\label{b3}
\left \{
\begin{array}{l}
\ep^2 u_{xx} + f(u,  (C(\ep) - \ep^2 u )/D ) = 0, \ \ \ x \in (0, 1)  \\[1ex]
u_x(0) = u_x(1) = 0, 
\end{array} 
\right.
\end{equation}
and 
\begin{equation}\label{b3_1}
\xi = \dis\frac{C(\ep)}{D} + \left( 1 - \frac{\ep^2}{D} \right) \int_0^1 u(x ; \ep)dx,
\end{equation}
respectively. Though the problem \eqref{b3} is quite similar to the one that was discussed in \cite{I}, there are two differences between them. First, the nonlinear term $f$ in \eqref{b3} does not depend on $x$ explicitly, and second the extra condition  \eqref{b3_1} is added. We assume that a solution $u$ of \eqref{b3} exhibits a sharp jump-up transition layer at $x= x^*(\ep) \in (0,1)$.\par 
To solve this problem, in Section~\ref{S2.1}, we assume that $C(\ep)$ and $ x^*(\ep)$ 
are represented by 
\begin{equation}\label{b8}
C(\ep) = C_0 + \ep C_1, \ (C_0/D \in I) 
\end{equation}
and
\begin{equation}\label{b9}
x^*(\ep) = x_0 +  \ep x_1,
\end{equation}
respectively
for small $\ep > 0$, where 
$I = (\un{v}, \ov{v})$ is given in the assumption (A1).
We note that $C_0$ and $x_0$ are to be determined by (A1), (A2) and \eqref{a6} in \eqref{b35x} and \eqref{b48_0}, respectively.
Since we consider an approximate solution of \eqref{b3} 
up to $O(\ep)$, the error terms $O(\ep^2)$ do not appear in the
right hand side of \eqref{b8} and \eqref{b9}, i.e., 
$C_1$ and $x_1$ are to be determined as functions of $\ep$, respectively, in Section~\ref{S2.4}.
We divide the interval $[0, 1]$ into two subintervals $[0, x^*(\ep)]$ and $[x^*(\ep), 1]$, and consider the following two boundary value problems: 
\begin{equation}\label{b5}
\left \{
\begin{array}{l}
\ep^2 u_{xx} + f( u,  (C(\ep) - \ep^2 u )/D ) = 0, \ \ \ x \in  (0, x^*(\ep))   \\[1ex]
u_x(0) = 0, \ u(x^*(\ep)) = \alpha
\end{array} 
\right.
\end{equation}
and
\begin{equation}\label{b6}
\left \{
\begin{array}{l}
\ep^2 u_{xx} + f( u,  (C(\ep) - \ep^2 u )/D ) = 0, \ \ \ x \in  (x^*(\ep), 1) \\[1ex]
u(x^*(\ep)) = \alpha, \ u_x(1) = 0, 
\end{array} 
\right.
\end{equation}
where $\alpha$ is an arbitrary fixed constant satisfying $h^-(C_0/D) < \alpha < h^+(C_0/D)$ and $C_0/D \in I$.
By using the singular perturbation method used in \cite{I}, we show the existence of solutions satisfying \eqref{b5} and \eqref{b6}. In Section \ref{S2.2}, we match these solutions in $C^1$-sense at $x = x^*(\ep)$, from which we find an approximate solution of \eqref{b3} up to $O(\ep)$ by determining the relations between $C_j$ and $x_j$ for $j=0,1$. Similarly, we use the equation \eqref{b3_1}, and obtain the other relations  between $C_j$ and $x_j$ for $j=0,1$ in Section \ref{S2.3}. 
Finally, in Section \ref{S2.4}, using the result in Sections \ref{S2.2} and \ref{S2.3}, we determine the unknown constants $C(\ep)$ and 
$x^*(\ep)$ uniquely, and obtain the desired result about the existence of single transition layer solutions. 

Here, we again note the difference between our problem and that in \cite{I}. In \cite{I}, since the nonlinear term $f$ does depend on $x$ explicitly, $x^*(\ep)$ is determined by the $C^1$-matching at $x = x^*(\ep)$ and then the condition \eqref{b3_1} is not necessary. While in our problem \eqref{b3}, since $f$ does not depend on $x$ explicitly and contains an unknown constant $C(\ep)$, then to determine $C(\ep)$ and $x^*(\ep)$ uniquely we need two relations;  the $C^1$-matching at $x = x^*(\ep)$ and the condition \eqref{b3_1}. This difference comes from the property of mass conservation. 

We use the following function spaces with positive $\ep$ in this section: 
\begin{equation*}
\begin{array}{lcl}
C^2_\ep[0,1] & = & \left\{ u \in C^2[0,1] \ | \ \dis \sum_{k=0}^2 \max_{0 \leq x \leq 1} \left| \left( \ep \frac {d}{dx}\right)^{\! k} \! u(x) \right| < \infty \right\}, \\[3.5ex]
\pot{C}^2_\ep[0,1] & = & \left\{ u \in C^2_\ep[0,1] \ | \ u_x(0) = 0, \ u_x(1) = 0 \right\},  \\[2ex]
C^2_{\ep,0}[0,1] & = & \left\{ u \in C^2_\ep[0,1] \ | \ u_x(0) = 0, \ u(1) = 0 \right\},  \\[2ex]
C^2_{\ep,1}[0,1] & = & \left\{ u \in C^2_\ep[0,1] \ | \ u(0) = 0, \ u_x(1) = 0 \right\}.
\end{array} 
\end{equation*}

\subsection{Solutions of \eqref{b5} and \eqref{b6}}\label{S2.1} 

First, we consider the approximation of the solution of \eqref{b5} up to $O(\ep)$. Applying the change of variables $x = x^*(\ep)y$, we have
\begin{equation}\label{b7}
\left \{
\begin{array}{l}
\ep^2 u_{yy} + (x^*(\ep))^2 f( u,  (C(\ep) - \ep^2 u )/D ) = 0, \ \ \ y \in (0, 1)   \\[1ex]
u_y(0) = 0, \ u(1) = \alpha.
\end{array} 
\right.
\end{equation}
In order to construct the outer approximation of the solution of \eqref{b7},
substituting 
\[
u(y) = U^-_0(y) + \ep U^-_1(y) + O(\ep^2)
\]
into \eqref{b7}, and comparing each coefficients of powers of $\ep$, we have 
\begin{equation}\label{b11}
f( U^-_0, C_0/D) = 0 
\end{equation}
and
\begin{equation}\label{b12}
f_u^- U^-_1 + f_v^- C_1/D = 0,
\end{equation}
where $f_u^- = f_u( U^-_0, C_0/D)$, $f_v^- = f_v( U^-_0, C_0/D)$. 
Since we consider a jump-up solution at $y=1$, it follows from \eqref{b11} and
\eqref{b12} that
\begin{equation}\label{b14}
U^-_0(y) = h^{-}(C_0/D)
\end{equation}
and
\begin{equation}\label{b15}
U^-_1(y;C_1) = - (f_v^- / f_u^-) \cdot (C_1/D) = h^-_v(C_0/D) \cdot (C_1/D),
\end{equation}
where we used the relation $f_u(h^-(v), v) h^-_v(v) + f_v(h^-(v), v) = 0$
obtained by the differentiation of $f( h^-(v), v) = 0$ in $v$.
It should be noted that $U^-_j$ $(j=0, 1)$ are constants independent of $y$. 
Since these outer approximations do not satisfy the boundary condition at $y = 1$, 
we must consider the correction of the above approximation to the solution of \eqref{b7} in a neighborhood of $y = 1$ 
with the aid of the inner approximation given by
\begin{equation}\label{b17}
\begin{array}{l}
u(y)  =  U^-_0(y) + \ep U^-_1(y;C_1) + \phi^-_0( \frac{y-1}{\ep}) 
 + \ep \phi^-_1(  \frac{y-1}{\ep} ) + O(\ep^2).
\end{array}
\end{equation}
Introducing the stretched coordinate $z = (y-1)/\ep$, and substituting 
\eqref{b17} into \eqref{b7}, and 
comparing each coefficients of powers of $\ep$, we have 
\begin{equation}\label{b18}
\left \{
\begin{array}{l}
\ddot{\phi}_0^- + x_0^2 \tilde{f}^- = 0, \ \ \ z \in (-\infty, 0)  
 \\[1ex]
\phi^-_0(-\infty) = 0, \ \phi^-_0(0) = \alpha - U^-_0(1) 
\end{array} 
\right.
\end{equation}
and
\begin{equation}\label{b19}
\left \{
\begin{array}{l}
\ddot{\phi}_1^- + x_0^2 \tilde{f}^-_u \phi^-_1 = F_1^-(z;C_1,x_1), \ \ \ z \in (-\infty, 0)  
 \\[1ex]
\phi^-_1(-\infty) = 0, \ \phi^-_1(0) = -U^-_1(1;C_1),
\end{array} 
\right.
\end{equation}
where the dot notation denotes $d/dz$, and 
$$
F_1^-(z;C_1,x_1) =  -2x_0 x_1 \tilde{f}^- - x_0^2 \tilde{f}^-_u U^-_1(1;C_1) -x_0^2 \tilde{f}^-_v C_1/D ,
$$
$$ \tilde{f}^- = f( h^-(C_0/D) + \phi^-_0, C_0/D ), \ 
\tilde{f}^-_u = f_u( h^-(C_0/D) + \phi^-_0, C_0/D ),
$$
and $\tilde{f}^-_v$ is similarly defined. 
From the assumptions (A1) and (A2), 
we find that \eqref{b18} has a unique monotone
increasing solution $\phi^-_0(z)$ since $\alpha > h^-(C_0/D)$ and $C_0/D \in I$. 
Moreover, we see that the solution of \eqref{b19} 
is explicitly given by
\begin{equation}\label{b21}
\begin{array}{l}
\phi^-_1(z;C_1,x_1) = -U^-_1(1;C_1) \dis\frac{ \dot{\phi}_0^-(z) }{ \dot{\phi}_0^-(0) } \\[1ex]
\hspace*{2cm} - \ \dot{\phi}_0^-(z)  \dis\int_z^0 \frac{1}{ (\dot{\phi}_0^-(\eta))^2} 
\dis\int_{-\infty}^{\eta}  \dot{\phi}_0^-(\zeta) 
F_1(\zeta;C_1,x_1) d\zeta d\eta. 
\end{array} 
\end{equation}
\par
Now, we put 
\begin{equation}\label{b22zz}
\begin{array}{lcl}
U^-(y;\ep;C_1,x_1) & = & U^-_0(y) + \ep U^-_1(y;C_1) \\[0.2cm]
 &  & + \ \theta(y)\phi^-_0( \frac{y-1}{\ep}) 
+  \ep \theta(y)\phi^-_1(  \frac{y-1}{\ep};C_1,x_1) ,
\end{array}
\end{equation}
where $\theta(y) \in C^{\infty}[0,1]$ is a cut-off function satisfying 
\begin{equation*}
\begin{array}{c}
\theta(y) = 0, \quad y \in [0, 1/2]; \quad \theta(y) = 1, \quad y \in [3/4, 1]; \\[1ex]
 0 \leqq \theta(y) \leqq 1, \quad y \in (1/2, 3/4). 
\end{array}
\end{equation*}
Since $U^-(y;\ep;C_1,x_1)$ is an $O(\ep)$ approximation to a solution of \eqref{b7}, we set
\begin{equation}\label{b22_0}
\tilde{u}^-(y;\ep;C_1,x_1) = U^-(y;\ep;C_1,x_1) + \ep r^-(y;\ep;C_1,x_1) 
\end{equation}
and rewrite \eqref{b7} as the following form with respect to the remainder term $r^-$: 
\begin{equation}\label{b22_1}
\left \{
\begin{array}{c}
\ep^2 r^-_{yy} + (x^*(\ep))^2 f(U^- +\ep r^- ,  (C(\ep) - \ep^2 U^- -\ep^3 r^- )/D ) /\ep \\[1ex]
+ \ \ep U^-_{yy}(y;\ep;C_1,x_1) \ = \ 0, \ \ \ y \in (0, 1)   \\[1ex]
r^-_y(0) = 0, \ r^-(1) = 0.
\end{array} 
\right.
\end{equation}
When we simply write \eqref{b22_1} as 
\begin{equation}\label{b22_2}
\begin{array}{c}
T(r^-;\ep;C_1,x_1) \ = \ 0, 
\end{array} 
\end{equation}
$T$ is a smooth mapping from $C^2_{\ep,0}[0,1]$ to $C[0,1]$, and then we have the following lemma: 
\begin{lemma}\label{lem1b} 
For any given constants $C_1^*$ and $x_1^*$, there exist $\ep_0 > 0$, $\rho_0 > 0$, and $K >0$ such that for any 
$\ep \in (0, \ep_0)$ and $(C_1,x_1) \in \Delta_{\rho_0} \equiv \{(C_1,x_1) \in {\bf R}^2 \ | \ |(C_1,x_1) - (C^*_1,x^*_1)| \leq \rho_0 \}$, 
\begin{description}
\item[(i)] \  $||T(0;\ep;C_1,x_1) ||_{C[0,1]} = o(1)$ uniformly in $(C_1,x_1) \in \Delta_{\rho_0}$ as $\ep \to 0$; 
\item[(ii)] \ for any $r_1, r_2 \in C^2_{\ep,0}[0,1]$, \\[1ex]
$\displaystyle{ \left|\left| \frac {\partial T}{\partial r^-}(r_1;\ep;C_1,x_1) -  \frac {\partial T}{\partial r^-}(r_2;\ep;C_1,x_1) \right|\right|_{C^2_{\ep,0}[0,1] \to C[0,1]}
 \ \leq K || r_1 - r_2||_{C^2_{\ep,0}[0,1]} }$; 
\item[(iii)] \ \hspace*{0.5cm}$ \displaystyle{ \left|\left|  \left( \frac {\partial T}{\partial r^-} \right)^{ \! \! -1}  \! \! (0;\ep;C_1,x_1) \right|\right|_{C[0,1] \to C^2_{\ep,0}[0,1]}  \leq K.  }$ 
\end{description}
Moreover, the results (i)-(iii) hold also for $\partial T/\partial C_1$ and $\partial T/\partial x_1$ in place of $T$.
\end{lemma}

Since this lemma is proved by the argument similar to that of \cite[Lemma 4.3]{MTH}, we omit the proof. 

Owing to Lemma \ref{lem1b}, we can apply the implicit function theorem to \eqref{b22_2}, and thus obtain the following: 
\begin{prop}\label{prop1b} 
There exist $\ep_1 > 0$ and $\rho_1 > 0$ such that for any $\ep \in (0,\ep_1)$ and $\rho \in \Delta_{\rho_1}$, there exists $r^-(y;\ep;C_1,x_1) \in C^2_{\ep,0}[0,1]$ satisfying 
$$   
T(r^-(y;\ep;C_1,x_1);\ep;C_1,x_1) \ = \ 0. 
$$
Moreover, $r^-(y;\ep;C_1,x_1)$, $\partial r^- \! / \partial C_1(y;\ep;C_1,x_1)$ and $\partial r^- \! / \partial x_1(y;\ep;C_1,x_1)$ are uniformly continuous with respect to $(\ep,C_1,x_1) \in (0, \ep_1) \times \Delta_{\rho_1}$ in $C^2_{\ep,0}[0,1]$-topology and satisfy 
\begin{equation*}
\left. 
\begin{array}{l}
 ||r^-(y;\ep;C_1,x_1) ||_{C^2_{\ep,0}[0,1]} \\[2ex]
\displaystyle{ \left|\left| \frac {\partial r^-}{\partial C_1}(y;\ep;C_1,x_1) \right|\right|_{C^2_{\ep,0}[0,1] }} \\[3ex]
\displaystyle{ \left|\left| \frac {\partial r^-}{\partial x_1}(y;\ep;C_1,x_1) \right|\right|_{C^2_{\ep,0}[0,1] }}
\end{array} 
\right\}
= o(1) \  \mbox{uniformly in} \ (C_1,x_1) \in \Delta_{\rho_1} \ \mbox{as} \ \ep \to 0. 
\end{equation*}
\end{prop}

Therefore, we obtain the solution of \eqref{b5} on $[0, x^*(\ep)]$, which takes the form
\begin{equation}\label{b34}
\begin{array}{l}
u^-(x ; \ep;C_1,x_1)  \ :=  \tilde{u}^-(\frac{x}{x^*(\ep)}; \ep; C_1,x_1) \\[1ex]
\hspace*{1cm}  = U^-_0( \frac{x}{x^*(\ep)}) + \ep U^-_1(\frac{x}{x^*(\ep)};C_1) + \theta(\frac{x}{x^*(\ep)}) \phi^-_0( \frac{x-x^*(\ep)}{\ep x^*(\ep)}) \\[1ex]
\hspace{1.2cm} 
+ \  \ep \theta(\frac{x}{x^*(\ep)}) \phi^-_1(   \frac{x-x^*(\ep)}{\ep x^*(\ep)};C_1,x_1) + \ep r^-(\frac{x}{x^*(\ep)};\ep;C_1,x_1 ), \\[1ex]
\hspace*{6cm} \ x \in [0, x^*(\ep)]. 
\end{array}
\end{equation}

Next, we consider the solution of \eqref{b6}.
Applying the change of variables $x = x^*(\ep) + (1-x^*(\ep))y$, we have
\begin{equation}\label{b23}
\left \{
\begin{array}{l}
\ep^2 u_{yy} + (1- x^*(\ep))^2 f( u,  (C(\ep) - \ep^2 u )/D ) = 0, \ \ \ y \in (0, 1)   \\[1ex]
u(0) = \alpha, \ u_y(1) = 0, \ 
\end{array} 
\right.
\end{equation}
where $C(\ep)$ and $x^*(\ep)$ are given by \eqref{b8} and \eqref{b9}, respectively.
Applying the same lines of argument as applied to \eqref{b7}, 
we can obtain the outer approximation of \eqref{b23} 
\[
u(y)  =  U^+_0(y) + \ep U^+_1(y;C_1) + O(\ep^2),
\]
where
\begin{equation}\label{b25}
U^+_0(y) = h^{+}(C_0/D) 
\end{equation}
and
\begin{equation}\label{b26}
U^+_1(y;C_1) = - (f_v^+ / f_u^+) \cdot (C_1/D) = h^+_v(C_0/D) \cdot (C_1/D),
\end{equation}
and $f_u^+ = f_u( U^+_0, C_0/D)$, $f_v^+ = f_v( U^+_0, C_0/D)$.
We note that $U^+_j$ $(j=0, 1)$ are constants independent of $y$.

Similarly to the problem \eqref{b7}, since these outer approximations do not satisfy the boundary condition at $y = 0$, 
we must consider the correction of the above approximation in a neighborhood of $y = 0$ 
with the aid of the inner approximation given by 
\begin{equation}\label{b28}
\begin{array}{l}
u(y)  =  U^+_0(y) + \ep U^+_1(y;C_1) + \phi^+_0( \frac{y}{\ep}) 
+  \ep \phi^+_1(  \frac{y}{\ep} ) + O(\ep^2).
\end{array}
\end{equation}
Introducing the stretched coordinate $z = y/\ep$, and substituting 
\eqref{b28} into \eqref{b23}, and 
comparing each coefficients of powers of $\ep$, we have 
\begin{equation}\label{b29}
\left \{
\begin{array}{l}
\ddot{\phi}_0^+ + (1-x_0)^2 \tilde{f}^+ = 0, \ \ \ z \in (0, \infty)  
 \\[1ex]
\phi^+_0(0) = \alpha - U^+_0(0), \  \phi^+_0(\infty) = 0 \ 
\end{array} 
\right.
\end{equation}
and
\begin{equation}\label{b30}
\left \{
\begin{array}{l}
\ddot{\phi}_1^+ + (1-x_0)^2 \tilde{f}^+_u \phi^+_1 = F_1^+(z;C_1,x_1), \ \ \ z \in (0, \infty)  
 \\[1ex]
 \phi^+_1(0) = -U^+_1(0;C_1), \ \phi^+_1(\infty) = 0, \
\end{array} 
\right.
\end{equation}
where the dot notation denotes $d/dz$, and 
$$
F_1^+(z;C_1,x_1) =  2(1-x_0) x_1 \tilde{f}^+ - (1- x_0)^2 \tilde{f}^+_u U^+_1(0;C_1) 
-(1-x_0)^2 \tilde{f}^+_v C_1/D ,
$$
$$ \tilde{f}^+ = f( h^+(C_0/D) + \phi^+_0, C_0/D ), \ 
\tilde{f}^+_u = f_u( h^+(C_0/D) + \phi^+_0, C_0/D ),
$$
and $\tilde{f}^+_v$ is similarly defined.
From the assumptions (A1) and (A2), 
we find that \eqref{b29} has a unique monotone
increasing solution $\phi^+_0(z)$ since $\alpha < h^+(C_0/D)$ and $C_0/D \in I$. 
Moreover, we see that the solution of \eqref{b30} is explicitly given by
\begin{equation}\label{b32}
\begin{array}{l}
\phi^+_1(z;C_1,x_1) = - \ U^+_1(0;C_1) \dis\frac{ \dot{\phi}_0^+(z) }{ \dot{\phi}_0^+(0) } \\[2ex]
\hspace*{2cm} - \ \dot{\phi}_0^+(z)  \dis\int_0^z \frac{1}{ (\dot{\phi}_0^+(\eta))^2} 
\dis\int_{\eta}^{\infty}  \dot{\phi}_0^+(\zeta) 
F_1^+(\zeta;C_1,x_1) d\zeta d\eta. 
\end{array} 
\end{equation}

Applying an argument similar to the above, we can find the solution 
$\tilde{u}^+(y;\ep;C_1,x_1)$ of \eqref{b23} as follows:
$$  
\tilde{u}^+(y;\ep;C_1,x_1)  =  U^+(y;\ep;C_1,x_1) + \ep r^+(y;\ep;C_1,x_1), 
$$
\begin{equation}\label{any1}
\begin{array}{l}
U^+(y;\ep;C_1,x_1) =  U^+_0(y) + \ep U^+_1(y;C_1)  \hspace*{1cm} \\[0.2cm]
\hspace*{1cm}+ \ \theta(1-y)\phi^+_0( \frac{y}{\ep}) 
+  \ep \theta(1-y)\phi^+_1(  \frac{y}{\ep};C_1,x_1 ). 
\end{array}
\end{equation}
Here, $r^+(y;\ep;C_1,x_1) \in C^2_{\ep,1}[0,1]$ satisfies 
\begin{equation*}
\left. 
\begin{array}{l}
 ||r^+(y;\ep;C_1,x_1) ||_{C^2_{\ep,1}[0,1]} \\[2ex]
\displaystyle{ \left|\left| \frac {\partial r^+}{\partial C_1}(y;\ep;C_1,x_1) \right|\right|_{C^2_{\ep,1}[0,1] }} \\[3ex]
\displaystyle{ \left|\left| \frac {\partial r^+}{\partial x_1}(y;\ep;C_1,x_1) \right|\right|_{C^2_{\ep,1}[0,1] }}
\end{array} 
\right\}
= o(1) \  \mbox{uniformly in} \ (C_1,x_1) \in \Delta_{\rho_2} \ \mbox{as} \ \ep \to 0. 
\end{equation*}
Moreover, 
$r^+(y;\ep;C_1,x_1)$, $\partial r^+ \! /\partial C_1(y;\ep;C_1,x_1)$ and $\partial r^+ \! /\partial x_1(y;\ep;C_1,x_1)$ are uniformly continuous with respect to $(\ep,C_1,x_1) \in (0, \ep_2) \times \Delta_{\rho_2}$
in $C^2_{\ep,1}[0,1]$-topology, 
where $\ep_2$ and $\rho_2$ are positive constants. 
Thus, we obtain the solutions of \eqref{b6} on $[x^*(\ep), 1]$ which takes the form 
\begin{equation}\label{b35}
\begin{array}{l}
\! \! \! u^+(x ; \ep;C_1,x_1)  \ :=  \tilde{u}^+(\frac{x-x^*(\ep)}{1- x^*(\ep)}; \ep;C_1,x_1) \\[1ex]
\hspace*{1cm}  = 
U^+_0( \frac{x-x^*(\ep)}{1- x^*(\ep)}) + \ep U^+_1(\frac{x-x^*(\ep)}{1- x^*(\ep)};C_1) + \theta(\frac{1-x}{1- x^*(\ep)}) \phi^+_0( \frac{x-x^*(\ep)}{\ep(1- x^*(\ep))}) 
 \\[1ex]
\hspace{1.2cm} 
+ \ \ep \theta(\frac{1-x}{1- x^*(\ep)}) \phi^+_1(  \frac{x-x^*(\ep)}{\ep(1- x^*(\ep))};C_1,x_1) + \ep r^+(y;\ep,C_1,x_1), \\[1ex]
\hspace*{6cm} x \in [x^*(\ep), 1]. 
\end{array}
\end{equation}

\subsection{$C^1$-matching at $x = x^*(\ep)$}\label{S2.2} 

We now consider the $C^1$-matching of $u^-(x ; \ep;C_1,x_1)$ and $u^+(x ; \ep;C_1,x_1)$ at $x = x^*(\ep)$ to obtain the approximation of the solution of \eqref{b3} up to $O(\ep)$. Since these two solutions are already continuous at $x = x^*(\ep)$, we then determine the values of $C_j$ and $x_j$ $(j=0, 1)$ in such a way that
\[
\Phi(\ep) := \ep x^*(\ep) ( 1 - x^*(\ep)) \{
\frac{d}{dx}u^-(x^*(\ep) ; \ep;C_1,x_1) - \frac{d}{dx}u^+(x^*(\ep) ; \ep;C_1,x_1) \} = o(\ep) 
\]
holds for small $\ep >0$. Noting that $U^-_j$ and $U^+_j$ are constants,  
it follows from \eqref{b34} and \eqref{b35} that
\[
\begin{array}{lcl}
\Phi(\ep) & = & \ep  ( 1 - x^*(\ep)) \{
\dot{\phi}^-_0(0)/\ep + \dot{\phi}^-_1(0;C_1,x_1) + o(1) \} \\[1ex]  
 & & - \ \ep  x^*(\ep) \{
\dot{\phi}^+_0(0)/\ep + \dot{\phi}^+_1(0;C_1,x_1)  + o(1) \} \\[1ex] 
& = & \Phi_0 + \ep \Phi_1(C_1,x_1) + o(\ep),
\end{array}
\]
where
\begin{equation}\label{b35zz}
\Phi_0 = (1-x_0) \dot{\phi}^-_0(0) - x_0  \dot{\phi}^+_0(0), 
\end{equation}
$$
\Phi_1(C_1,x_1) =  (1-x_0) \dot{\phi}^-_1(0;C_1,x_1) - x_1  \dot{\phi}^-_0(0) - x_0  \dot{\phi}^+_1(0;C_1,x_1)
- x_1  \dot{\phi}^+_0(0).
$$

First, we consider $\Phi_0 =  (1-x_0) \dot{\phi}^-_0(0) - x_0  \dot{\phi}^+_0(0) = 0 $.
It follows from \eqref{b18} that 
\[
\begin{array}{rcl}
0 & = & \dis\int_{-\infty}^0 \{ \ddot{\phi}^-_0 \dot{\phi}^-_0
+ x_0^2 f( h^-(C_0/D) + \phi^-_0, C_0/D ) \dot{\phi}^-_0 \} dz 
\\[3ex]
& = & 
\dis\frac{(\dot{\phi}^-_0(0))^2}{2} + x_0^2
\dis\int_{h^-(C_0/D)}^\alpha f(u, C_0/D ) du,
\end{array}
\]
which implies 
\[
\dot{\phi}^-_0(0) = x_0 \sqrt{-2 \dis\int_{h^-(C_0/D)}^\alpha f(u, C_0/D ) du}. 
\]
Similarly, it follows from \eqref{b29} that 
\[
\dot{\phi}^+_0(0) = (1-x_0) \sqrt{ 2 \dis\int_\alpha^{h^+(C_0/D)} f(u, C_0/D ) du}. 
\]
Therefore, we have
\[
\Phi_0 = 
-\dis\frac{ 2 x_0(1-x_0) J(C_0/D)}
{ \sqrt{- 2 \dis\int_{h^-(C_0/D)}^\alpha f(u, C_0/D ) du}
+ \sqrt{  2 \dis\int_\alpha^{h^+(C_0/D)} f(u, C_0/D ) du} },
\]
where $J=J(v)$ is given by \eqref{a4}. Hence, noting the assumptions (A.1) and (A.2),
it follows from \eqref{b14}, \eqref{b25} and $\Phi_0 = 0$ that
\begin{equation}\label{b35x}
C_0 = Dv^* \ \ \ \text{and} \ \ \ U^{\pm}_0(y) = h^{\pm}(v^*).
\end{equation}
We note that though $\Phi_0$ depends on both $C_0$ and $x_0$, the solution satisfying $\Phi_0 = 0$ is determined by only $C_0 = D v^*$ for any $x_0 \in (0,1)$.
Moreover, we have
\begin{equation}\label{b36}
\dot{\phi}^-_0(0) = x_0 \sqrt{-2 \dis\int_{h^-(v^*)}^\alpha f(u, v^* ) du}
\end{equation}
and
\begin{equation}\label{b37}
\dot{\phi}^+_0(0) = (1-x_0) \sqrt{ 2 \dis\int_\alpha^{h^+(v^*)} f(u, v^* ) du}. 
\end{equation}

Next, we consider
$\Phi_1(C_1,x_1) = (1-x_0) \dot{\phi}^-_1(0;C_1,x_1) - x_1  \dot{\phi}^-_0(0) - x_0  \dot{\phi}^+_1(0;C_1,x_1) - x_1  \dot{\phi}^+_0(0) = 0$.
Note the following relation: 
\[
\dis\int_{-\infty}^0 x_0^2 \tilde{f}^-_u  \dot{\phi}^-_0 dz
=
 - \dis\int_{-\infty}^0 \dddot{\phi}^-_0  dz
= 
-  \ddot{\phi}^-_0 (0), 
\]
where we used the relation $\dddot{\phi}^-_0 + x_0^2 \tilde{f}_u^- \dot{\phi}_0^- = 0$ obtained by
the differentiation of the first equation of \eqref{b18} in $z$. 
Then, it follows from \eqref{b21} that
\begin{equation}\label{b38}
\begin{array}{l}
 \dot{\phi}^-_1(0;C_1,x_1) 
= \dis  -U^-_1(1;C_1)\frac {\ddot{\phi}^-_0 (0)}{\dot{\phi}_0^-(0)} + 
 \frac{1}{\dot{\phi}_0^-(0)}\dis\int_{-\infty}^0 F_1^-(z;C_1,x_1) \dot{\phi}^-_0 dz
\\[3ex]
\ \ =
\dis\frac{1}{\dot{\phi}_0^-(0)} \left(   -U^-_1(1;C_1)\ddot{\phi}^-_0 (0) 
-2x_0 x_1 \dis\int_{-\infty}^0  \tilde{f}^- \dot{\phi}^-_0 dz \right. \\[3ex]
\left. \hspace*{2.5cm} - x_0^2 U^-_1(1;C_1) \dis\int_{-\infty}^0  \tilde{f}^-_u \dot{\phi}^-_0 dz
 - \frac{ x_0^2 C_1 }{D} \dis\int_{-\infty}^0 \tilde{f}^-_v \dot{\phi}^-_0 dz \right)
\\[3ex]
\ \ = 
\dis\frac{1}{\dot{\phi}_0^-(0)} \left(
-2x_0 x_1 \dis\int_{h^-(v^*)}^\alpha  f(u, v^*) du
- \dis\frac{ x_0^2 C_1 }{D} \dis\int_{h^-(v^*)}^\alpha  f_v(u, v^*) du \right). 
\end{array}
\end{equation}
Similarly, it follows from \eqref{b29} and \eqref{b32} that
\begin{equation}\label{b39}
\begin{array}{l}
\dot{\phi}^+_1(0;C_1,x_1)  
= 
\dis\frac{1}{\dot{\phi}_0^+(0)} \left(
-2(1-x_0 ) x_1 \dis\int^{h^+(v^*)}_\alpha  f(u, v^*) du \right.
\\[3ex]
\ \ \ \ \ \ \ \ \ \ \ \ \ \ \ \ \ \ \ \ \ \ \ \ \ \ \ \ \ \ \ \ \ \ \ 
\left. + \dis\frac{ (1-x_0)^2 C_1 }{D} \dis\int^{h^+(v^*)}_\alpha  f_v(u, v^*) du \right).
\end{array}
\end{equation}
Hence, by \eqref{b38} and \eqref{b39}, we have
\[
\begin{array}{rcl}
\Phi_1(C_1,x_1) & = & (1-x_0) \dot{\phi}^-_1(0;C_1,x_1) - x_1  \dot{\phi}^-_0(0) - x_0  \dot{\phi}^+_1(0;C_1,x_1) - x_1  \dot{\phi}^+_0(0)
\\[2ex]
& =: & K(x_0)x_1 + M(x_0)C_1 + R(x_0),
\end{array}
\]
where
\begin{equation}\label{b40}
\begin{array}{l}
K(x_0) = \dis\frac{x_0(1-x_0)}{\dot{\phi}^-_0(0)}
\left( -2 \dis\int_{h^-(v^*)}^\alpha f(u, v^* ) du \right) - \dot{\phi}^-_0(0)
\\[3ex]
\ \ \ \ \ \ \ \ \ \ \ 
+ \dis\frac{x_0(1-x_0)}{\dot{\phi}^+_0(0)}
\left( 2 \dis\int_\alpha^{h^+(v^*)} f(u, v^* ) du \right) - \dot{\phi}^+_0(0),
\end{array}
\end{equation}
\begin{equation}\label{b41}
\begin{array}{l}
M(x_0) = - \dis\frac{x_0^2(1-x_0)}{D\dot{\phi}^-_0(0)}
\left(  \dis\int_{h^-(v^*)}^\alpha f_v(u, v^* ) du \right) 
\\[3ex]
\ \ \ \  \ \ \ \ \ \ \ \ \ \ 
- \dis\frac{x_0(1-x_0)^2}{D\dot{\phi}^+_0(0)}
\left( \dis\int_\alpha^{h^+(v^*)} f_v(u, v^* ) du \right),
\end{array}
\end{equation}
and $R(x_0) = 0$. Moreover, it follows from \eqref{b36}, \eqref{b37} and 
\eqref{b40} that
\[
\begin{array}{rcl}
K(x_0) & = & (1-2x_0) \left( \sqrt{ -2 \dis\int_{h^-(v^*)}^\alpha f(u, v^* ) du } 
- \sqrt{ 2 \dis\int_\alpha^{h^+(v^*)} f(u, v^* ) du } \right)
\\[4ex]
& = & \dis\frac{ -2(1-2x_0) J(v^*)}{ \sqrt{ -2 \dis\int_{h^-(v^*)}^\alpha f(u, v^* ) du }
+ \sqrt{ 2 \dis\int_\alpha^{h^+(v^*)} f(u, v^* ) du } } = 0,
\end{array}
\]
where $J(v)$ is given by \eqref{a4}. 
This implies that $\Phi_1$ does not depend on $x_1$.
Noting
\[
\dis\int_\alpha^{h^+(v^*)} f(u, v^* ) du = - \dis\int_{h^-(v^*)}^\alpha f(u, v^* ) du
\]
by $J(v^*) = 0$, 
it follows from \eqref{b36}, \eqref{b37} and 
\eqref{b41} that
\[
M(x_0)  = 
- \dis\frac{x_0(1-x_0)}{D  \sqrt{ -2 \dis\int_{h^-(v^*)}^\alpha f(u, v^* ) du } }
\dis\int_{h^-(v^*)}^{h^+(v^*)} f_v(u, v^* ) du \neq 0 
\]
by \eqref{a4x}. Therefore, it follows from \eqref{b15}, \eqref{b26} and $\Phi_1 = 0$
that 
\begin{equation}\label{b43x}
C_1 = -R(x_0)/M(x_0) = 0 \ \ \ \ \text{and} \ \ \ \  U^{\pm}_1(y;C_1) = h_v^{\pm}(v^*)C_1/D = 0.
\end{equation}

\subsection{Computation of \eqref{b3_1}}\label{S2.3} 

To complete the construction of the approximate solution of \eqref{b3} satisfying \eqref{b3_1}, 
we determine the values of $C_j$ and $x_j$ $(j=0,1)$ by the conservation law \eqref{b3_1}. 
Here, we put 
\begin{equation}\label{b44}
\Psi(\ep) := \dis\frac{C(\ep)}{D} + \left( 1 - \frac{\ep^2}{D} \right) \int_0^1 u(x;\ep;C_1,x_1)dx - \xi. 
\end{equation}
\eqref{b3_1} is equivalent to $\Psi(\ep) = 0$. 
Using \eqref{b34}, \eqref{b35}, \eqref{b35x} and \eqref{b43x}, we 
have
\begin{equation}\label{b45}
\int_0^1 u(x;\ep;C_1,x_1)dx  = \int_0^{x^*(\ep)} u^-(x;\ep;C_1,x_1)dx  
+ \int_{x^*(\ep)}^1 u^+(x;\ep;C_1,x_1)dx,
\end{equation}
where
\[
\begin{array}{l}
\dis\int_0^{x^*(\ep)} u^-(x;\ep;C_1,x_1)dx  = x^*(\ep) \{
\dis\int_0^1 ( U^-_0 +\ep  U^-_1 + o(\ep) ) dy 
\\[3ex]
\ \ \ \ \ \ \ \ \ \ \ \ \ \ \ 
+ \, \ep \dis\int_{-\infty}^0  ( {\phi}_0^-(z) +  \ep {\phi}_1^-(z;C_1,x_1) + o(\ep) )dz
\, \} 
\\[3ex]
\ \ \ \ \ 
= ( x_0 + \ep x_1 + o(\ep) )\{ h^-(v^*) + \ep h_v^-(v^*)C_1/D 
+ \ep \dis\int_{-\infty}^0   {\phi}_0^-(z) dz + o(\ep) \}
\end{array}
\]
and
\[
\begin{array}{l}
\dis\int_{x^*(\ep)}^1 u^+(x;\ep;C_1,x_1)dx  = (1-x^*(\ep)) \{
\dis\int_0^1 ( U^+_0 +\ep  U^+_1 + o(\ep) ) dy 
\\[3ex]
\ \ \ \ \ \ \ \ \ \ \ \ \ \ \ 
+ \, \ep \dis\int^{\infty}_0  ( {\phi}_0^+(z) +  \ep {\phi}_1^+(z;C_1,x_1) + o(\ep))dz
\, \} 
\\[3ex]
\ \ \ \ \ 
= (1- x_0 - \ep x_1 + o(\ep) )\{ h^+(v^*) + \ep h_v^+(v^*)C_1/D 
+  \ep \dis\int^{\infty}_0   {\phi}_0^+(z) dz + o(\ep) \}.
\end{array}
\]
Substituting \eqref{b8} and \eqref{b45} into \eqref{b44}, we have  
\begin{equation}\label{b45_1}
\begin{array}{l}
\Psi(\ep) = \dis \{v^* + x_0 h^-(v^*) + (1- x_0) h^+(v^*) - \xi \} \\[1ex]
+ \ \ep \left\{ x_0 h_v^-(v^*)C_1/D + 
x_0 \dis\int_{-\infty}^0   {\phi}_0^-(z) dz + x_1 h^-(v^*)\right. \\[3ex]
\left. + \ (1-x_0) h_v^+(v^*)C_1/D 
+ (1-x_0) \dis\int^{\infty}_0   {\phi}_0^+(z) dz   - x_1 h^+(v^*) \right\} 
+ \ o(\ep) \\[3ex]
 =: \Psi_0 + \ep \Psi_1(C_1,x_1) + o(\ep). 
\end{array}
\end{equation}
Comparing each coefficients of powers of $\ep$ in \eqref{b45_1}, we have 
$ \Psi_i = 0 \ (i=0,1)$. 
Noting \eqref{a6} and $h^-(v^*) < h^+(v^*)$ by the assumption (A1), 
it follows from $\Psi_0 = 0$ that
\begin{equation}\label{b48_0}
\begin{array}{l}
x_0 = x^* = \dis\frac{ h^+(v^*) + v^* - \xi}{h^+(v^*) - h^-(v^*) }
\ \ \ \text{and} \ \ \ 
0 < x_0 < 1.
\end{array}
\end{equation}
Moreover, it follows from $\Psi_1(C_1,x_1) = 0$ that
\begin{equation}\label{b48_1}
\begin{array}{l}
\{x_0 h_v^-(v^*)+(1-x_0)h_v^+(v^*)\}C_1/D \hspace{1cm} \\[1ex]
\hspace{2cm} -  \ \{h^+(v^*)-h^-(v^*)\}x_1 + I_1(x_0) = 0,
\end{array}
\end{equation}
where 
\[
I_1(x_0) =  x_0 \dis\int_{-\infty}^0   {\phi}_0^-(z) dz 
+ (1-x_0) \dis\int^{\infty}_0   {\phi}_0^+(z) dz
\]
is a function of $x_0$. Since $C_1 = 0$ by \eqref{b43x}, 
it follows from \eqref{b48_1} that 
$x_1$ is uniquely determined by 
\begin{equation}\label{b48_2}
\begin{array}{l}
 x_1 = \dis \frac {I_1(x_0)}{h^+(v^*)-h^-(v^*)}.
\end{array}
\end{equation}
Thus, we see that \eqref{b34} and \eqref{b35} give 
the approximation of the solution of \eqref{b3} with \eqref{b3_1}. 

\subsection{Determination of $C(\ep)$ and $x^*(\ep)$}\label{S2.4} 
Finally, we determine  $C(\ep)$ and $x^*(\ep)$ uniquely such that \eqref{b3} with \eqref{b3_1} have a 
single transition layer solution $u(x;\ep)$ at the layer position $ x = x^*(\ep)$.

First, the coefficients $C_i$ and $x_i$ $(i=0,1)$ are determined step by step as follows:  $C_0$ (by \eqref{b35x}) $\longrightarrow$ 
$x_0$ (by \eqref{b48_0})  $\longrightarrow$ $C_1$ (by \eqref{b43x})  $\longrightarrow$ $x_1$ (by \eqref{b48_2}). We note that $\Phi(\ep) = o(\ep)$ and $\Psi(\ep) = o(\ep)$; $\Phi$ and $\Psi$ are not identically zero for these $C_1$ and $x_1$. 

Next, we set $(C_1^*, x_1^*) = (C_1, x_1)$ in Lemma \ref{lem1b}, and 
consider $C(\ep) = C_0 + \ep \bar{C}_1$ and 
$x(\ep) = x_0 + \ep \bar{x}_1$. 
We can take $(\bar{C}_1, \bar{x}_1)$ around $(C_1^*, x_1^*)$ 
so as to satisfy $\Phi(\ep) = 0$ and $\Psi(\ep) = 0$ as follows:
Let us define $\Phi^*(\bar{C}_1, \bar{x}_1;\ep)$ and  
$\Psi^*(\bar{C}_1, \bar{x}_1;\ep)$ by $\Phi(\ep) = \ep \Phi^*(\bar{C}_1, \bar{x}_1;\ep)$ and $\Psi(\ep) = \ep \Psi^*(\bar{C}_1, \bar{x}_1;\ep)$, respectively. We easily find that there exist two positive constants $\delta$ and $\ep_3 (< \min \{\ep_1, \ep_2\})$ 
such that $\Phi^*(\bar{C}_1, \bar{x}_1;\ep)$ and $\Psi^*(\bar{C}_1, \bar{x}_1;\ep)$ are continuous in $\bar{C}_1 \in (C_1^*-\delta, C_1^*+\delta), \bar{x}_1 \in (x_1^*-\delta, x_1^*+\delta)$ and $\ep \in [0,\ep_3)$, and are $C^1$-functions of $\bar{C}_1$ and $\bar{x}_1$. Moreover, we can easily find that 
\[ \left\{ 
\begin{array}{l}
\Phi^*(C_1^*, x_1^*;0) = 0, \
\dis\frac {\partial \Phi^*}{\partial \bar{C}_1}(C_1^*, x_1^*;0) = M(x_0), \
\dis\frac {\partial \Phi^*}{\partial \bar{x}_1}(C_1^*, x_1^*;0) = 0, 
\\[3ex]
\dis\Psi^*(C_1^*, x_1^*;0) = 0, \ 
\dis\frac {\partial \Psi^*}{\partial \bar{C}_1}(C_1^*, x_1^*;0) 
= ( x_0 h_v^-(v^*) +  (1 - x_0) h_v^+(v^*))/D,  
\\[3ex]
\dis\frac {\partial \Psi^*}{\partial \bar{x}_1}(C_1^*, x_1^*;0) = h^+(v^*) - h^-(v^*), 
\end{array} \right.
\]
which implies that 
\[
\begin{array}{l}
\dis\frac {\partial (\Phi^*, \Psi^*)}{\partial (\bar{C}_1, \bar{x}_1)}(C_1^*, x_1^*;0) 
 = M(x_0) (h^+(v^*) - h^-(v^*)) \  \ne \ 0. 
\end{array}
\]
Then, we can apply the implicit function theorem to $\Phi^*(\bar{C}_1, \bar{x}_1;\ep) = 0$ and  
$\Psi^*(\bar{C}_1, \bar{x}_1;\ep) = 0$, and find that there exist $\bar{C}_1 = \bar{C}_1(\ep)$
and $\bar{x}_1 = 
\bar{x}_1(\ep)$ for $\ep \in [0,\ep_3)$ satisfying $\bar{C}_1(0) = C_1^*, \ \bar{x}_1(0) = 
 x_1^*$,  
$\Phi^*(\bar{C}_1, \bar{x}_1;\ep) = 0$, and $\Psi^*(\bar{C}_1, \bar{x}_1;\ep) = 0$.

Substituting $C_1 = \bar{C}_1(\ep)$ and $x_1 = \bar{x}_1(\ep)$ into 
\eqref{b8}, \eqref{b9}, \eqref{b34} and \eqref{b35}, we obtain the following existence result: 

\begin{theo}\label{th2}
Suppose the assumptions (A1) and (A2).
For any $\ep \in (0, \ep_3)$, there exists a family of solutions $u^\ep(x) \in \pot{C}^2_\ep[0,1]$ of \eqref{b3} and \eqref{b3_1}. Furthermore, the following estimate holds:
\begin{equation*}
\begin{array}{l}
\dis \left|\left| u^\ep(x) - U^-(\frac x{x^*(\ep)};\ep;  \bar{C}_1(\ep), \bar{x}_1(\ep)) \right|\right|_
{C^2_{\ep,0}[0,x^*(\ep)]} \\[4ex]
+ \ \dis \left|\left| u^\ep(x) - U^+(\frac {x-x^*(\ep)}{1-x^*(\ep)};\ep; \bar{C}_1(\ep), \bar{x}_1(\ep)) \right|\right|_{C^2_{\ep,1}[x^*(\ep),1]} = o(\ep)
\end{array}
\end{equation*}
as $\ep \to 0$. 
\end{theo}

\begin{remark}
In Section \ref{S2.1}, we took $U_0^-(y) = h^-(C_0/D)$ in \eqref{b14} and $U_0^+(y) = h^+(C_0/D)$ in \eqref{b25} as a jump-up solution. If we take $U_0^-(y) = h^+(C_0/D)$ in \eqref{b14} and $U_0^+(y) = h^-(C_0/D)$ in \eqref{b25}, we can obtain a jump-down solution $u^\ep(x)$ at $x = x^*(\ep)$. 
\end{remark}

\begin{remark}
We can perform the matched asymptotic expansion of any order $k$, i.e., we can obtain 
a family of solutions $u^\ep(x) \in \pot{C}^2_\ep[0,1]$ of \eqref{b3} and \eqref{b3_1}
with the estimate
\begin{equation*}
\begin{array}{l}
\dis \left|\left| u^\ep(x) - U^-_k(\frac x{x^*(\ep)};\ep;  \bar{C}_k(\ep), \bar{x}_k(\ep)) \right|\right|_
{C^2_{\ep,0}[0,x^*(\ep)]} \\[4ex]
+ \ \dis \left|\left| u^\ep(x) - U^+_k(\frac {x-x^*(\ep)}{1-x^*(\ep)};\ep; \bar{C}_k(\ep), \bar{x}_k(\ep)) \right|\right|_{C^2_{\ep,1}[x^*(\ep),1]} = o(\ep^k)
\end{array}
\end{equation*}
as $\ep \to 0$, where $U^{-}_k$ and $U^+_k$ are the $k$th order approximate solutions 
given by equations similar
to \eqref{b22zz} and \eqref{any1}, respectively.
\end{remark}

By using the relation $v^\ep(x) = (C(\ep) - \ep^2 u^\ep(x))/D$, we obtain 

\begin{cor}\label{cor1} \rm
For sufficiently small $\ep >0$, \eqref{a1} has a family of  stationary solutions $(u^\ep(x), v^\ep(x))$ 
satisfying \eqref{a2xyz} with the following properties:

\no (i)
\begin{equation}\label{b101}
\dis\lim_{\ep \to 0} v^\ep(x) = V^*(x) 
\end{equation}
uniformly on $x \in [0, 1]$, where $V^*(x)$ is given by \eqref{a9}.

\no (ii)
For each $\sigma > 0$ with $ 0 < \sigma < \min(x^*, 1-x^*)$, 
\begin{equation}\label{b102}
\dis\lim_{\ep \to 0} u^\ep(x) = U^*(x) 
\end{equation}
uniformly on $x \in [0, x^* - \sigma] \cup [x^* + \sigma, 1]$,
where $U^*(x)$ is given by \eqref{a8}.

\no (iii) Let 
$\tilde{u}^\ep(\zeta) = u^\ep(x^*(\ep)+\ep \zeta)$ and 
$\tilde{v}^\ep(\zeta) = v^\ep(x^*(\ep)+\ep \zeta)$, where
$\zeta = (x - x^*(\ep))/\ep$ is the stretched coordinate around the $C^1$-matching point 
$x = x^*(\ep)$. Then, for each $\sigma > 0$, 
\[
\dis\lim_{\ep \to 0} \tilde{v}^\ep(\zeta) = v^*
\]
and
\[
\dis\lim_{\ep \to 0} \tilde{u}^\ep(\zeta) = Q(\zeta)
\]
uniformly on $\zeta \in [- \sigma, \sigma]$, where $Q(\zeta)$ satisfies
\begin{equation}\label{b103y}
\left \{
\begin{array}{l}
Q_{\zeta \zeta} + f( Q, v^*) = 0, \\[1ex]
Q(\pm \infty) = h^{\pm}(v^*), \ \ Q(0) = h^0(v^*). 
\end{array} 
\right.
\end{equation}
\end{cor}

\begin{remark}\label{vco} \rm
Corollary~\ref{cor1}(ii) shows that the $u$-component of the stationary solution
exhibits a jump-up transition with $O(1)$-amplitude
when $\ep$ is sufficiently small. In contrast,
the $v$-component of the stationary solution
exhibits a jump-down transition with $O(\ep^2)$-amplitude
because $v^\ep(x) = ( C(\ep) - \ep^2 u^\ep(x))/D$ by \eqref{b2zz}, where
$C(\ep)$ defined by \eqref{b8} is constant. 
The matched asymptotic expansion of the second order 
can give the coefficient of $O(\ep^2)$-term
of $C(\ep)$, which is lengthy and involved in general. 
However, in Section~3, we present a simple example such that the $O(\ep^2)$-term 
of the $v$-component of a stationary solution can be clearly expressed.
\end{remark}

{\bf Proof of Corollary~\ref{cor1}}. (i) Since $u^\ep(x)$ is bounded, it follows from \eqref{b8} and \eqref{b35x} that
$\lim_{\ep \to 0} v^\ep(x) = \lim_{\ep \to 0} (C(\ep) - \ep^2 u^\ep(x))/D = v^*$
uniformly on $x \in [0, 1]$. 

\no
(ii) It follows from \eqref{b9}, \eqref{b22zz}, \eqref{b35x} and \eqref{b48_0} that
$$
\lim_{\ep \to 0} U^-(\frac x{x^*(\ep)};\ep;  \bar{C}_1(\ep), \bar{x}_1(\ep)) = h^-(v^*)
$$
uniformly on $x \in [0, x^* - \sigma]$. Similarly, we see that
$$
\lim_{\ep \to 0} U^+(\frac {x-x^*(\ep)}{1-x^*(\ep)};\ep; \bar{C}_1(\ep), \bar{x}_1(\ep)) 
= h^+(v^*)
$$
uniformly on $x \in [x^* + \sigma, 1]$. Therefore, we have
\eqref{b102} by Theorem~\ref{th2}.

\no
(iii) 
It is clear that 
$\lim_{\ep \to 0} \tilde{v}^\ep(\zeta) = v^*$
uniformly on $\zeta \in [- \sigma, \sigma]$ because of the assertion (i).
When $\zeta \leq 0$, noting $x = x^*(\ep) + \ep \zeta$, $x= x^*(\ep)y$,
$z = (y-1)/\ep$ and \eqref{b35x}, we have
$$
\begin{array}{l}
\tilde{u}^\ep(\zeta)
 =  u^-( x^*(\ep) + \ep \zeta ; \ep, \bar{C}_1(\ep), \bar{x}_1(\ep))
= \tilde{u}^-(1 + \ep z ; \ep; \bar{C}_1(\ep), \bar{x}_1(\ep)) 
\\[1ex]
\ \ \ \ \ \ \ \, =  U^-( 1 + \ep z ; \ep, \bar{C}_1(\ep), \bar{x}_1(\ep)) + O(\ep)
= U_0^-(1) + \phi_0^-(z) + O(\ep) 
\\[1ex]
\ \ \ \ \ \ \ \, =  h^-(v^*) + \phi_0^-(\frac{\zeta}{x_0}) + O(\ep),
\end{array} 
$$ 
where $u^-$, $\tilde{u}^-$, $U^-$ and $\phi_0^-$ are given by \eqref{b34}, \eqref{b22_0},
\eqref{b22zz} and \eqref{b18}, respectively. On the other hand, applying the change of 
variable $z = \eta/x_0$ to \eqref{b18}, we find that 
$Q^-(\eta) := h^-(v^*) + \phi_0^-(\eta/x_0)$ satisfies 
\[
\left \{
\begin{array}{l}
Q^-_{\eta \eta} + f( Q^-, v^*) = 0, \\[1ex]
Q^-(- \infty) = h^{-}(v^*), \ \ Q^-(0) = h^0(v^*), \ \ 
Q^-_{\eta}(0) = \dot{\phi}^-_0(0)/x_0, 
\end{array} 
\right.
\]
where we choose $\alpha = h^0(v^*) \in (h^-(v^*), h^+(v^*))$. 
Similarly, when $\zeta \geq 0$, we have 
$\tilde{u}^\ep(\zeta) = h^+(v^*) + \phi_0^+(\zeta/(1-x_0)) + O(\ep)$, and 
find that $Q^+(\eta) := h^+(v^*) + \phi_0^+(\eta/(1-x_0))$
satisfies 
\[
\left \{
\begin{array}{l}
Q^+_{\eta \eta} + f( Q^+, v^*) = 0, \\[1ex]
Q^+(+ \infty) = h^{+}(v^*), \ \ Q^+(0) = h^0(v^*), \ \ 
Q^+_{\eta}(0) = \dot{\phi}^+_0(0)/(1-x_0).
\end{array} 
\right.
\]
We can obtain 
$$
Q(\eta) = 
\left \{
\begin{array}{l}
Q^-(\eta)  \ \ \ (\eta \leq 0) \\[1ex]
Q^+(\eta) \ \ \ (\eta > 0) 
\end{array} 
\right.
$$ 
satisfying the same equation as
\eqref{b103y} by the $C^1$-matching of 
$Q^-(\eta)$ and $Q^+(\eta)$ at $\eta = 0$ because 
$\dot{\phi}^-_0(0)/x_0  = \dot{\phi}^+_0(0) /(1-x_0)$ by $\Phi_0 = 0$, where
$\Phi_0$ is given by \eqref{b35zz}. 
Since we can identify $\zeta$ with $\eta$
in the limit of $\ep \to 0$, by using a similar argument in the proof of \cite[Lemma 1.1]{NF},
we see that 
$\lim_{\ep \to 0} \tilde{u}^\ep(\zeta) = Q(\zeta)$ 
uniformly on $\zeta \in [- \sigma, \sigma]$, where $Q(\zeta)$ satisfies \eqref{b103y}.
\Qed

\section{Stability analysis of transition layer solutions} \label{stability} 

In this section, we perform the stability analysis of $(u^\ep(x), v^\ep(x))$ given by
Corollary~\ref{cor1} under a constraint derived from 
the conservation law \eqref{a2}.

We consider the linearized eigenvalue problem
\begin{equation}\label{c1}
\mathcal{L}^\ep \Phi^\ep = \lambda^\ep \Phi^\ep, \ \ \
\Phi^\ep = 
\left(
\begin{array}{c}
\varphi^\ep \\
\psi^\ep
\end{array}
\right) \in X, 
\ \ \
\mathcal{L}^\ep =
\left(
\begin{array}{cc}
L^\ep & f_v^\ep \\
- f_u^\ep & M^\ep
\end{array}
\right),
\end{equation}
under the Neumann boundary condition, where 
$f_u^\ep = f_u ( u^\ep(x), v^\ep(x) )$, 
$f_v^\ep = f_v ( u^\ep(x), v^\ep(x) )$, 
\[
L^\ep := \ep^2 \dis\frac{d^2}{dx^2} + f_u ( u^\ep(x), v^\ep(x) ), \ \ 
M^\ep := D \dis\frac{d^2}{dx^2} - f_v ( u^\ep(x), v^\ep(x) ),
\]
\begin{equation}\label{c2}
X = \{ \ (\varphi, \psi) \in L^2(0, 1) \times L^2(0, 1) \ | \ \int_0^1 (\varphi + \psi) dx = 0 \ \},
\end{equation}
and the domain $D(\mathcal{L}^\ep)$ is naturally defined. 
We note that the constrained condition in the definition $X$ is naturally 
derived from the conservation law \eqref{a2}.

In this section, we use the notations
\begin{equation}\label{c2xx}
f_u^*(x) = f_u( U^*(x), V^*(x)) \ \ 
\text{and} \ \ 
f_v^*(x) = f_v( U^*(x), V^*(x)),
\end{equation}
where $(U^*(x), V^*(x))$ is given by \eqref{a8} and \eqref{a9}. 
From the assumptions (A1) and (A3), we have
\begin{equation}\label{fa1}
f_u^*(x) < 0 \ \  \ \ ( 0 \leq x \leq 1)
\end{equation}
and
\begin{equation}\label{fa2}
f_v^*(x) - f_u^*(x) > 0 \ \  \ \ ( 0 \leq x \leq 1),
\end{equation}
respectively. Moreover, we denote the inner product in 
$L^2(0,1)$ by $\lg \, \cdot \, , \, \cdot \, \rg$.

\subsection{Preliminaries}   

We denote by $\mu_j^\ep$ $(j =0,1,2 \cdots )$ the eigenvalues of 
$L^\ep : L^2(0, 1) \to L^2(0, 1)$ under the Neumann boundary condition;
\[
\sigma(L^\ep) = \{ \mu_j^\ep \}_{j=0}^\infty,  \ \ \ 
\mu_0^\ep > \mu_1^\ep > \cdots > \mu_j^\ep \to - \infty \quad (j \to \infty).
\]

We note that 
$V^*(x) \equiv v^*$ satisfies $V^*_x(x^*) = 0$,
where $V^*(x)$ is given by \eqref{a9} and $V^*(x) = \lim_{\ep \to 0} v^\ep(x) $
by \eqref{b101}. 
Moreover, we note that $f(h^\pm(v^*), v^*) = 0$ by the assumption (A1).
Noting Corollary~\ref{cor1}, 
we see that the following propositions are true by
\cite[Proposition 5.1]{HS} and
\cite[Corollary 1.3, Lemma 1.4 and Lemma 2.3]{NF}. 

\begin{prop}\label{prop1c}
Under the assumptions (A1) and (A2), for sufficiently small $\ep > 0$, the following properties hold:

\no
{\rm (i)}  Let $\phi_0^\ep(x)$ be the $L^2$-normalized eigenfunction of $L^\ep$
corresponding to $\mu_0^\ep$. Then, 
\begin{equation}\label{c3}
\mu_0^\ep = O(e^{-C/\ep}) \ \  \text{as} \ \ \ep \to 0
\end{equation}
for some $C>0$, and
\begin{equation}\label{c4}
\lim_{\ep \to 0} \frac{1}{\sqrt{\ep}} \int_0^1 \phi_0^\ep(x) dx 
= \kappa^*(h^+(v^*) - h^-(v^*)),
\end{equation}
where 
\begin{equation}\label{c4x}
\kappa^* = \left( \int_{-\infty}^\infty Q_\zeta(\zeta)^2 d\zeta  \right)^{-1/2}
\end{equation}
and $Q=Q(\zeta)$ is the solution of \eqref{b103y}. 

\no
{\rm (ii)} 
For each $p \in H^1(0, 1)$, we have
\begin{equation}\label{c6}
\dis\lim_{\ep \to 0} \left\lg p, \left(\frac{\phi_0^\ep}{\sqrt{\ep}} \right) \! f_u^\ep \right\rg = 0
\end{equation}
and
\begin{equation}\label{c7}
\dis\lim_{\ep \to 0} \left\lg p, \left(\frac{\phi_0^\ep}{\sqrt{\ep}} \right) \! f_v^\ep \right\rg = 
 \kappa^* p(x^*)J'(v^*),
\end{equation}
where $\kappa^*$ and $J'(v^*)$ are given by \eqref{c4x}
and \eqref{a4x}, respectively.
\end{prop}

\begin{prop}\label{prop2c} 
Under the assumptions (A1) and (A2), for sufficiently small $\ep > 0$, the following properties hold:

\no
{\rm (i)}  There exist $\mu_* > 0$ independent of $\ep$
such that $ \mu_1^\ep <  -\mu_* < \mu_0^\ep$. Moreover,
$(L^\ep - \lambda )^{-1}P^\ep : L^2(0, 1) \to L^2(0, 1)$ is well-defined 
for any $\lambda \in \Lambda = 
\{ \, \lambda \in \C \  | \ \mbox{\rm Re} \lambda > - \mu_* \, \}$, and
\begin{equation}\label{c5z}
|| (L^\ep - \lambda )^{-1}P^\ep q ||_{L^2} \leq \dis\frac{1}{|\lambda - \underline{\mu}|} ||q||_{L^2} \ \  \text{for} \ q \in L^2(0,1)
\end{equation}
holds for any given $\underline{\mu} \in ( \mu_1^\ep, -\mu_*)$,
where $P^\ep$ is the orthonormal projection onto the orthogonal complement of $\phi_0^\ep$.

\no
{\rm (ii)} For $q \in L^\infty(0, 1)$ and $\lambda \in \Lambda$, 
\begin{equation}\label{c5}
\dis\lim_{\ep \to 0} (L^\ep - \lambda)^{-1} P^\ep q = \frac{ q(x)}{f_u^*(x) - \lambda}
\ \ \text{strongly in} \ L^2(0, 1).
\end{equation}
Moreover, the convergence in \eqref{c5}
is uniform with respect to $\lambda \in \Lambda$ and $q$ on a $H^1$-bounded set.
\end{prop}

\subsection{Properties of eigenvalues} \rm

Noting \eqref{fa1}, we define
\begin{equation}\label{c8}
\nu = \min\{ \mu_*  , \ \inf_{0 \leq x \leq 1} (- f_u^*(x) ) \} /2  > 0,
\end{equation}
where $\mu_*$ is given by Proposition~\ref{prop2c}(i), and 
$f_u^*(x)$ and $f_v^*(x)$ are given by \eqref{c2xx}. 
We consider the eigenvalues of \eqref{c1} in
$\Lambda_* = 
\{ \, \lambda \in \C \  | \ \mbox{\rm Re} \lambda > - \nu \, \}$.

Decomposing the first component of the eigenfuction of
$\mathcal{L}^\ep$ as $\varphi^\ep = a^\ep \phi_0^\ep + w^\ep$,
the eigenvalue problem \eqref{c1} is rewritten as

\begin{equation}\label{c9}
\left \{
\begin{array}{l}
a^\ep(\mu_0^\ep - \lambda^\ep) = - \lg \psi^\ep, f_v^\ep \phi^\ep_0 \rg  \\[2ex]
( L^\ep - \lambda^\ep ) w^\ep = - P^\ep(f_v^\ep \psi^\ep) \\[2ex]
( M^\ep - \lambda^\ep) \psi^\ep - f_u^\ep w^\ep = a^\ep f_u^\ep \phi^\ep_0,
\end{array} 
\right. 
\end{equation}

\no
where $(\mu_0^\ep, \phi_0^\ep(x))$ is the principal eigenpair of
$L^\ep$, $a^\ep \in \mbox{\bf C}$, $w^\ep$ satisfies $\lg w^\ep, 
\phi_0^\ep \rg =0$, and $P^\ep$ is the orthogonal projection
onto the orthogonal complement of $\phi_0^\ep$. 
We note that $(\varphi^\ep, \psi^\ep) = ( a^\ep \phi_0^\ep + w^\ep, \psi^\ep) \in X$, which implies 
that the constrained condition
\begin{equation}\label{c9x}
a^\ep \int_0^1 \phi_0^\ep dx + \int_0^1 w^\ep dx  +   \int_0^1 \psi^\ep dx = 0
\end{equation}
holds. This equation plays a key role to characterize a critical eigenvalue
which essentially determines the stability of the stationary solutions.
Noting $\lambda^\ep \in \Lambda_* \subset \Lambda$, it follows from Proposition~\ref{prop2c}(i) that

\begin{equation}\label{c10}
w^\ep = - ( L^\ep - \lambda^\ep )^{-1} P^\ep(f_v^\ep \psi^\ep).
\end{equation}

\no
Substituting \eqref{c10} into the third equation in \eqref{c9}, we have

\begin{equation}\label{c11}
( M^\ep - \lambda^\ep) \psi^\ep + f_u^\ep ( L^\ep - \lambda^\ep )^{-1} P^\ep(f_v^\ep \psi^\ep) = a^\ep f_u^\ep \phi^\ep_0.
\end{equation}

In what follows, we suppose $\lg \psi^\ep, \psi^\ep \rg = 1$. 
In fact, 
if $\psi^\ep \equiv 0$, then
$w^\ep \equiv 0$ and $a^\ep =0$ by \eqref{c4}, \eqref{c9x} and \eqref{c10}.

\begin{lemma}\label{lem1c}
Under the assumptions (A1) and (A2),
for each $p \in H^1(0, 1)$, we have
\begin{equation}\label{c12}
- \lg D \psi^\ep_{x}, p_x \rg + \lg - f_v^\ep \psi^\ep - \lambda^\ep \psi^\ep
 + \frac{f_u^\ep f_v^\ep }{ f_u^* - \lambda^\ep }\psi^\ep, p \rg = o(1)
\ \ \text{as} \ \ \ep \to 0.
\end{equation}
\end{lemma}

{\bf Proof}. 
Noting $|| \psi^\ep ||_{L^2} = 1$, we have
$$
\left| \int_0^1 \psi^\ep dx \right| \leq || \psi^\ep ||_{L^2} = 1.
$$
Moreover, noting $\lambda^\ep \in \Lambda_* \subset \Lambda$, it follows from \eqref{c5z} and
\eqref{c10} that
$$
 || w^\ep ||_{L^2}  
= || ( L^\ep - \lambda^\ep )^{-1} P^\ep(f_v^\ep \psi^\ep)||_{L^2} 
\leq \dis\frac{1}{|\lambda^\ep - \underline{\mu}|} ||f_v^\ep \psi^\ep||_{L^2}
$$
holds for any given $\underline{\mu} \in ( \mu_1^\ep, -\mu_*)$.
Since $\text{Re}\lambda^\ep > - \nu \geq -\mu_*/2$ by 
$\lambda^\ep \in \Lambda_*$, we have
$$
|| w^\ep ||_{L^2} \leq  \frac{2}{\mu_*} || f^\ep_v \psi^\ep ||_{L^2},
$$
which implies
$$
\left| \int_0^1 w^\ep dx \right| \leq  || w^\ep ||_{L^2} \leq \frac{2}{\mu_*} || f^\ep_v \psi^\ep ||_{L^2} 
\leq \frac{2}{\mu_*} || f^\ep_v ||_{L^\infty}.
$$
Since $||f^\ep_v||_{L^\infty}$ is uniformly bounded with respect to $\ep$, it follows from \eqref{c4} and \eqref{c9x} that
$$
\left| a^\ep \! \! \int_0^1 \phi_0^\ep dx \right| = \left| \int_0^1 w^\ep dx  +   \int_0^1 \psi^\ep dx \right| \leq 
1 + \frac{2}{\mu_*} || f^\ep_v ||_{L^\infty} 
$$
and $| a^\ep | \leq C/\sqrt{\ep}$ for sufficiently small $\ep$, where $C$ is a positive constant independent of $\ep$.
Thus, it follows from \eqref{c6} that 
$$
| \lg a^\ep f_u^\ep \phi^\ep_0, p \rg | = 
| a^\ep |  \sqrt{\ep} \, \Big{|} \Big{\lg} p, \left(\frac{\phi_0^\ep}{\sqrt{\ep}} \right) \! f_u^\ep \Big{\rg} \Big{|} = o(1) 
\ \ \text{as} \ \ \ep \to 0.
$$
Thus, by using integration by parts, we obtain \eqref{c12} due to \eqref{c5} and \eqref{c11}. \Qed

\begin{lemma}\label{lem2c}
Under the assumptions (A1) and (A2), there exists a positive constant
$R$ independent of $\ep$ such that
$|\lambda^\ep| \leq R$ holds.
\end{lemma}

{\bf Proof}. It follows from \eqref{c12} that
$$
- \lg D \psi^\ep_{x}, \psi^\ep_x \rg + \lg - f_v^\ep \psi^\ep - \lambda^\ep \psi^\ep
 + \frac{f_u^\ep f_v^\ep }{f_u^* - \lambda^\ep}\psi^\ep, \psi^\ep \rg = o(1)
\ \ \text{as} \ \ \ep \to 0.
$$
By using $\lg \psi^\ep, \psi^\ep \rg =1$, we have
$$
\lambda^\ep = - D \lg  \psi^\ep_x, \psi^\ep_x \rg
- \lg f_v^\ep \psi^\ep, \psi^\ep \rg 
+ \lg \dis\frac{f_u^\ep f_v^\ep }{f_u^* - \lambda^\ep}\psi^\ep, \psi^\ep \rg + o(1) \ \ \text{as} \ \ \ep \to 0.
$$
Hence, noting $\lambda^\ep \in \Lambda_*$,
$| \lg f_v^\ep \psi^\ep, \psi^\ep \rg | \leq  || f_v^\ep ||_{L^\infty} ||\psi^\ep ||_{L^2}^2$ and  $|| \psi^\ep ||_{L^2} = 1$, we have
$$
\begin{array}{l}
-\nu \ < \ \text{Re}\lambda^\ep = - D \lg  \psi^\ep_x, \psi^\ep_x \rg
- \lg f_v^\ep \psi^\ep, \psi^\ep \rg 
+ \text{Re} \lg \dis\frac{f_u^\ep f_v^\ep }{f_u^* - \lambda^\ep}\psi^\ep, \psi^\ep \rg + o(1) 
\\[3ex]
\ \ \ \ \ \ \ \ \ \ \ \ \ \ \ \ \ \,  \leq   \ || f_v^\ep ||_{L^\infty}  
+ \text{Re}\lg \dis\frac{f_u^\ep f_v^\ep }{f_u^* - \lambda^\ep}\psi^\ep, \psi^\ep \rg + o(1) \ \ \text{as} \ \ \ep \to 0
\end{array}
$$
and
$$
\text{Im}\lambda^\ep  =  \text{Im}\lg \dis\frac{f_u^\ep f_v^\ep }{f_u^* - \lambda^\ep}\psi^\ep, \psi^\ep \rg + o(1) \ \ \text{as} \ \ \ep \to 0.
$$
Here, we assume that $\lim_{\ep \to 0} |\lambda^\ep| = +\infty$. 
Since
$$
| f_u^* - \lambda^\ep| \geq |\lambda^\ep| - | f_u^*| \geq |\lambda^\ep| -
|| f_u^* ||_{L^\infty} > 0
$$
holds for sufficiently small $\ep$, noting $|| \psi^\ep ||_{L^2} = 1$, we have
$$
| \lg \dis\frac{f_u^\ep f_v^\ep }{f_u^* - \lambda^\ep}\psi^\ep, \psi^\ep \rg |
\leq \int_0^1 \dis\frac{|f_u^\ep | |f_v^\ep| }{|f_u^* - \lambda^\ep|} |\psi^\ep|^2 dx
\leq \dis\frac{ ||f_u^\ep ||_{L^\infty}  ||f_v^\ep||_{L^\infty}  }{ |\lambda^\ep| -
|| f_u^* ||_{L^\infty} } \to 0
$$
as $\ep \to 0$. This implies that $|\lambda^\ep|$ is bounded for sufficiently small 
$\ep$, which leads to a contradiction. Thus,
we see that $\lim_{\ep \to 0} |\lambda^\ep| < +\infty$ holds. \Qed \\

Under the assumptions (A1) and (A2), 
it follows from Lemma~\ref{lem2c} that there are no eigenvalues in $
\{ \, \lambda \in \C \ | \ \text{Re} \lambda > -\nu, \  |\lambda| > R \, \}$,
where $\nu$ is given by \eqref{c8}.
In what follows, we consider eigenvalues in 
$\Lambda_\delta = \{ \, \lambda \in \C \ | \ \text{Re} \lambda \geq -\delta, \ |\lambda| \leq R \, \}$, where $\delta$ 
is an arbitrarily given constant satisfying $0 < \delta \leq \nu$.
Noting that $\lambda^\ep$ is continuous in $\ep$ on a compact set 
$\Lambda_\delta$, we suppose that $\lim_{\ep \to 0} \lambda^\ep$ exists in $ \Lambda_\delta$ because we are interested in the behavior of eigenvalues which 
determine the stability of $(u^\ep(x), v^\ep(x))$.

\begin{theo}\label{th3} 
Under the assumptions (A1), (A2) and (A3), any eigenvalue 
$\lambda^\ep \in \Lambda_{\delta_1}$ must satisfy
$$ \dis\lim_{\ep \to 0} \text{\rm Re} \lambda^\ep \leq 0 $$
and 
$$ \dis\lim_{\ep \to 0} \text{\rm Re} \lambda^\ep < 0 \ \ 
\text{if} \ \ \dis\lim_{\ep \to 0} \text{\rm Im} \lambda^\ep \neq 0, $$ 
where $\delta_1 = \min\{ \nu, \dis\inf_{0 \leq x \leq 1 } \dis\frac{f_v^*(x) - f_u^*(x)}{2} \} > 0$.
\end{theo}

{\bf Proof}. 
We note that $\delta_1$ is well-defined by \eqref{fa2} and \eqref{c8}.
It follows from \eqref{c12} that
$$
- \lg D \psi^\ep_{x},  \psi^\ep_{x} \rg + \lg - f_v^\ep \psi^\ep - \lambda^\ep \psi^\ep
 + \frac{f_u^\ep f_v^\ep }{f_u^* - \lambda^\ep}\psi^\ep, \psi^\ep \rg = o(1)
\ \ \text{as} \ \ \ep \to 0,
$$
which implies
\[
D \lg  \psi^\ep_x, \psi^\ep_x \rg  + \lambda^\ep \Big{\lg} 
\frac{\lambda^\ep + (f_v^\ep - f_u^\ep )}{\lambda^\ep - f_u^*}\psi^\ep, \psi^\ep \Big{\rg} = o(1)
\ \ \text{as} \ \ \ep \to 0.
\]

Let $\alpha^\ep = \text{Re} \lambda^\ep$ and $\beta^\ep = \text{Im} \lambda^\ep$. Then, we have
$$
\begin{array}{l}
\Big{\lg} \dis\frac{\lambda^\ep + (f_v^\ep - f_u^\ep )}{\lambda^\ep - f_u^*}\psi^\ep, \psi^\ep \Big{\rg}
=
\Big{\lg} \dis\frac{ (f_v^\ep - f_u^\ep ) + \alpha^\ep + i \beta^\ep}{(\alpha^\ep - f_u^*) + i \beta^\ep}\psi^\ep, \psi^\ep \Big{\rg} \\[3ex]
\ \ \ \ 
=
\Big{\lg} \dis\frac{ (\alpha^\ep - f_u^* )^2 +  f^\ep_v (\alpha^\ep - f_u^* ) + (\beta^\ep)^2
- i  f^\ep_v \beta^\ep}{(\alpha^\ep - f_u^*)^2 + (\beta^\ep)^2}\psi^\ep, \psi^\ep \Big{\rg} + o(1)
\end{array}
$$
as $\ep \to 0$ because $\lim_{\ep \to 0} f_u^\ep(x) = f_u^*(x)$ uniformly on 
$x \in [0, x^* - \sigma] \cup [x^* + \sigma, 1]$ for any 
$\sigma > 0$ with $ 0 < \sigma < \min(x^*, 1-x^*)$ due to Corollary~\ref{cor1}.
Let 
$$
K^\ep_1 :=  (\alpha^\ep - f_u^* )^2 +  f^\ep_v (\alpha^\ep - f_u^* ) + (\beta^\ep)^2 \ \
\text{and} \ \  K^\ep_2 := (\alpha^\ep - f_u^*)^2 + (\beta^\ep)^2. 
$$
Since $\alpha^\ep > -\delta_1$ and \eqref{c8} by $\lambda^\ep \in \Lambda_{\delta_1}$, when $\ep$ is sufficiently small, it follows from \eqref{fa1} that 
$$
\alpha^\ep - f_u^* > -  f_u^*/3 \geq \delta_1/3
$$
holds, which implies 
$K^\ep_2 \geq \delta_1^2/9 > 0$.
Similarly, when $\ep$ is sufficiently small, it follows from \eqref{fa2} that 
$$
\begin{array}{l}
 \alpha^\ep + (f_v^\ep - f_u^* ) =  \alpha^\ep - \{ - (f_v^\ep - f_u^*) \} 
\\[2ex]
\ \ \ \ \ = 
 \alpha^\ep - \{ - (f_v^* - f_u^*) \} + (f_v^\ep - f_v^*)
\\[2ex]
\ \ \ \ \  > (f_v^* - f_u^* )/3 \geq \delta_1/3
\end{array}
$$
on $x \in [0, x^* - \sigma] \cup [x^* + \sigma, 1]$ 
for any $\sigma > 0$ with $ 0 < \sigma < \min(x^*, 1-x^*)$
because 
$\lim_{\ep \to 0} f_v^\ep(x) = f_v^*(x)$ uniformly on 
$x \in [0, x^* - \sigma] \cup [x^* + \sigma, 1]$  due to Corollary~\ref{cor1}.
Therefore, we have 
$$
K^\ep_1 = (\alpha^\ep - f_u^* ) \{ \alpha^\ep + (f_v^\ep - f_u^* ) \} + (\beta^\ep)^2 \geq \delta_1^2/9 > 0
$$ 
on $x \in [0, x^* - \sigma] \cup [x^* + \sigma, 1]$ for sufficiently small $\ep$.
In addition, by Lemma~\ref{lem2c} and $|| f_v^\ep ||_{L^\infty} = O(1)$ as
$\ep \to 0$, we have
\begin{equation}\label{c13x}
0 < C_1 < \Big{\lg} \frac{K_1^\ep }{K_2^\ep}\psi^\ep, \psi^\ep \Big{\rg} < C_2
\end{equation}
for sufficiently small $\ep$, where $C_1$ and $C_2$ are independent of $\ep$.
Then, we have
$$
\begin{array}{l}
\lambda^\ep \Big{\lg} 
\dis\frac{\lambda^\ep + (f_v^\ep - f_u^\ep )}{\lambda^\ep - f_u^*}\psi^\ep, \psi^\ep \Big{\rg} 
= 
(\alpha^\ep + i \beta^\ep ) \Big{\lg} 
\frac{K_1^\ep - i f_v^\ep \beta^\ep}{K_2^\ep}\psi^\ep, \psi^\ep \Big{\rg} + o(1)
\\[3ex]
\ \ \ \ 
= \alpha^\ep \Big{\lg} \dis\frac{K_1^\ep }{K_2^\ep}\psi^\ep, \psi^\ep \Big{\rg} 
+  (\beta^\ep)^2 \Big{\lg}  \frac{f_v^\ep }{K_2^\ep}\psi^\ep, \psi^\ep \Big{\rg} 
\\[3ex]
\ \ \ \ \ \ \ \ \ \ \ 
+ i \beta^\ep \Big{\{} - \alpha^\ep \Big{\lg}  \dis\frac{f_v^\ep }{K_2^\ep}\psi^\ep, \psi^\ep \Big{\rg} 
+ \Big{\lg} \frac{K_1^\ep }{K_2^\ep}\psi^\ep, \psi^\ep \Big{\rg} \Big{\}} + o(1) 
\ \ \text{as} \ \ \ep \to 0,
\end{array}
$$
which implies
\[
D \lg  \psi^\ep_x, \psi^\ep_x \rg + \alpha^\ep \Big{\lg} \dis\frac{K_1^\ep }{K_2^\ep}\psi^\ep, \psi^\ep \Big{\rg} 
+  (\beta^\ep)^2 \Big{\lg}  \frac{f_v^\ep }{K_2^\ep}\psi^\ep, \psi^\ep \Big{\rg}  = o(1)
\]
and
\[
\beta^\ep \Big{\{} - \alpha^\ep \Big{\lg}  \dis\frac{f_v^\ep }{K_2^\ep}\psi^\ep, \psi^\ep \Big{\rg} 
+ \Big{\lg} \frac{K_1^\ep }{K_2^\ep}\psi^\ep, \psi^\ep \Big{\rg} \Big{\}} = o(1) 
\]
as $\ep \to 0$. Therefore, we have
\begin{equation}\label{c14}
\alpha^\ep D \lg  \psi^\ep_x, \psi^\ep_x \rg 
+ \{ (\alpha^\ep)^2 +  (\beta^\ep)^2 \} \Big{\lg} \dis\frac{K_1^\ep }{K_2^\ep}\psi^\ep, \psi^\ep \Big{\rg} 
= o(1)
\end{equation}
as $\ep \to 0$. 
Noting \eqref{c13x}, it follows from \eqref{c14} that
$$
\left( \alpha^\ep + \frac{D \lg  \psi^\ep_x, \psi^\ep_x \rg }{2 \big{\lg} (K_1^\ep/K_2^\ep) \psi^\ep, \psi^\ep \big{\rg} } 
\right)^2 + ( \beta^\ep )^{2} -
\left( \frac{D \lg  \psi^\ep_x, \psi^\ep_x \rg }{2  \big{\lg} (K_1^\ep/K_2^\ep) \psi^\ep, \psi^\ep\big{\rg} } \right)^2 = o(1)
$$
as $\ep \to 0$. In the limit of $\ep \to 0$, this equation seemingly expresses
a circle included in the 
left half plane in $\mbox{\bf C}$, and the circle is tangential to the imaginary axis at the origin.
Thus, we obtain Theorem~\ref{th3}.  \Qed \\

\begin{theo}\label{th4}
Under the assumptions (A1), (A2) and (A3), the stationary solutions
$(u^\ep(x), v^\ep(x))$ are stable if $J'(v^*) > 0$. i.e., any eigenvalue $\lambda^\ep$
of the linearized operator of \eqref{a1} at $(u^\ep(x), v^\ep(x))$ in $X$ given by \eqref{c2}
must satisfy $\text{Re}\lambda^\ep < 0$ for sufficiently small $\ep  > 0$ if $J'(v^*) > 0$.
\end{theo}

{\bf Proof}.
By Theorem~\ref{th3}, it suffices to prove that 
any eigenvalue
$\lambda^\ep \in \Lambda_{\delta_1}$ satisfying
$\lim_{\ep \to 0} \lambda^\ep = 0$ must satisfy 
$\text{Re}\lambda^\ep < 0$ for sufficiently small $\ep  > 0$ if $J'(v^*) > 0$.
It follows from \eqref{c12} that
$$
\lg D \psi^*_{x}, p_x \rg = 0 \ \ \text{for} \  \forall p \in H^1(0,1),
$$
where $\psi^* = \lim_{\ep \to 0} \psi^\ep$.
Hence, noting $\lg \psi^\ep, \psi^\ep \rg = 1$, we have $\psi^* \equiv 1$. 
Moreover, by \eqref{c5} and \eqref{c10}, we have $w^* = - f^*_v/f^*_u$, where
$w^* = \lim_{\ep \to 0} w^\ep$. Therefore, it follows from 
\eqref{c4} and \eqref{c9x} that
$$
\begin{array}{rcl}
\dis\lim_{\ep \to 0} a^\ep \sqrt{\ep} & = & - \dis\lim_{\ep \to 0} 
\dis\frac{  \dis\int_0^1 w^\ep dx  +   \int_0^1 \psi^\ep dx }
{ \dis\frac{1}{\sqrt{\ep}} \int_0^1 \phi_0^\ep dx }
\\[5ex]
& = &
 \dis\frac{1 }
{ \kappa^* (h^+(v^*) - h^-(v^*))}  \int_0^1 \dis\frac{f_v^* - f_u^*}{f_u^*} dx.
\end{array}
$$
Thus, by \eqref{c3}, \eqref{c7} and the first equation of \eqref{c9}, we obtain
\begin{equation}\label{cri}
\lambda^\ep = \frac{1}{a^\ep} \lg \psi^\ep, f_v^\ep \phi^\ep_0 \rg + \mu_0^\ep = 
 \ep (\kappa^*)^2 \cdot \dis\frac{ h^+(v^*) - h^-(v^*) }{  \dis\int_0^1 \dis\frac{f_v^* - f_u^*}{f_u^*} dx } 
\cdot J'(v^*) + o(\ep)
\end{equation}
as $\ep \to 0$. 
Thus, by \eqref{fa1} and \eqref{fa2}, we see that
$\lambda^\ep$ satisfying
$\lim_{\ep \to 0} \lambda^\ep = 0$
must satisfy $\mbox{\rm Re} \lambda^\ep < 0$ for 
sufficiently small $\ep$ if $J'(v^*) > 0$.
\Qed

\begin{remark}\label{rem3y}
Notice that \eqref{cri} gives the characterization of an eigenvalue satisfying 
$\lim_{\ep \to 0} \lambda^\ep = 0$. This eigenvalue is called the {\it critical eigenvalue}
which essentially determines the stability of $(u^\ep(x), v^\ep(x))$.
Although our result suggests that the stationary solutions $(u^\ep(x), v^\ep(x))$ are unstable if $J'(v^*) < 0$, we cannot conclude that it is true
because we have not yet proved the existence of the critical eigenvalue in $\Lambda_{\delta_1}$.
To prove this existence, 
we must consider the solvability of \eqref{c11} 
with respect to $\psi^\ep$ 
when $\lambda^\ep = O(\ep)$ as $\ep \to 0$.
In the case of reaction-diffusion systems of FitzHugh-Nagumo type studied 
by \cite{HS,NF}, such solvability problems can be solved by 
the Lax-Milgram theorem \cite{NS} which is the most powerful and
standard tool for 
solving the linear elliptic PDEs.  
In contrast to their cases, to solve
our problem by the Lax-Milgram theorem, we need 
the higher order terms of the convergence rates 
of \eqref{c6} and \eqref{c5} in $\ep$, respectively.
These terms are to be investigated for future studies.

\end{remark}

\section{Examples} 

In this section, we present some helpful examples for understanding our results;
we consider \eqref{a1} on an interval $0 < x < 1$ under the Neumann boundary condition,
where $0 < \ep \ll D$ and $f(u,v)$ is a smooth function with bistable nonlinearity. In our numerical simulations,
the values of the diffusion coefficients $\ep$ and $D$ are given by $\ep = 0.01$ and $D=1.0$, respectively. 

\subsection{Stable case}

We present mass-conserving reaction-diffusion 
systems with bistable nonlinearity, which have
stable stationary solutions with a single internal transition layer.
We note that their stability is guaranteed by Theorem~\ref{th4}.
The profiles of the stationary solutions in the figures in this subsection  
are numerically obtained by applying the Newton method to \eqref{b1};
the initial values of the iterations are 
given by the final states of numerical solutions of the time
evolutional PDE \eqref{a1}.

\begin{figure}[ht!]
\centering
\includegraphics[width=12cm, bb=0 0 908 830]{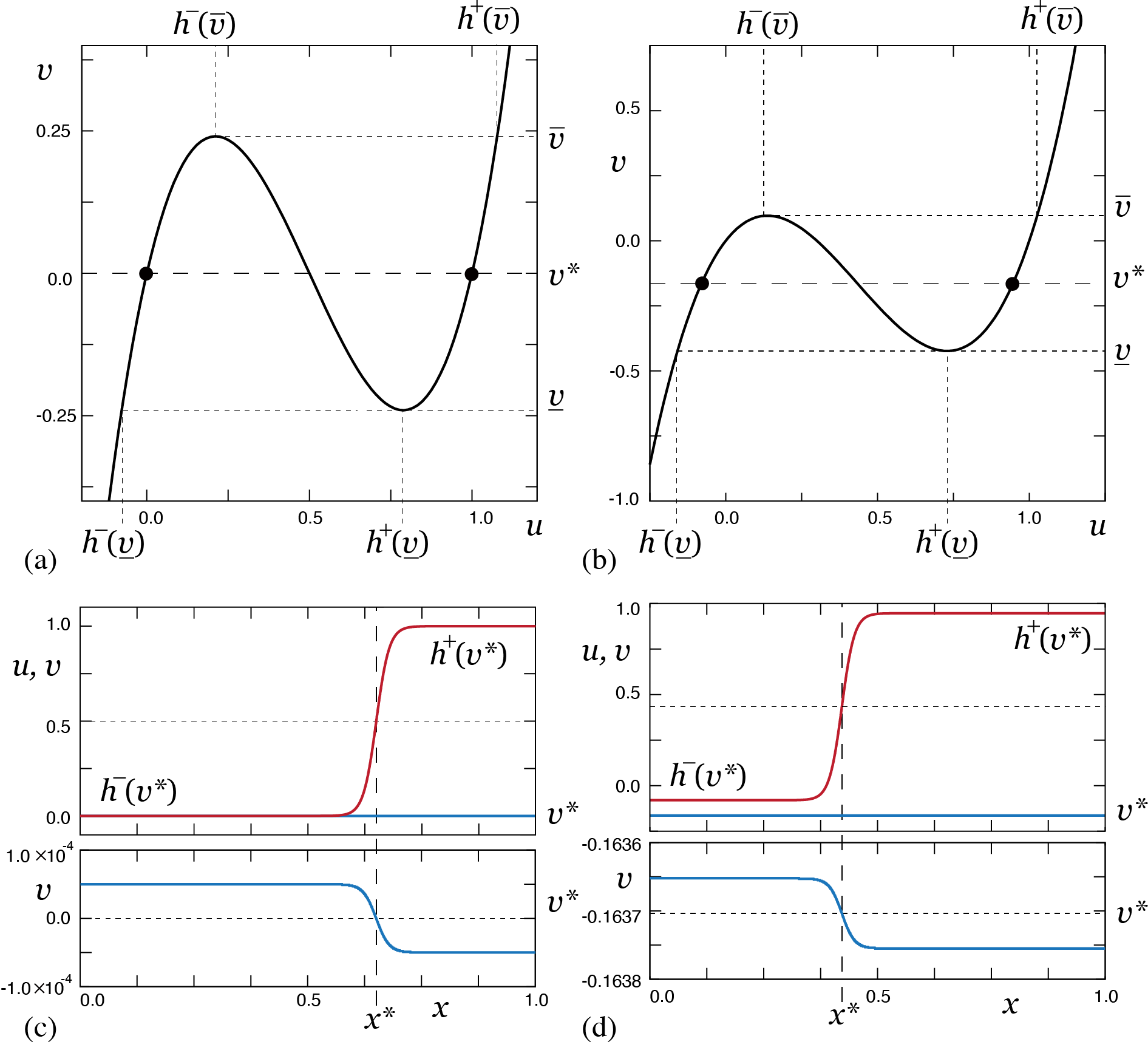}
  \caption{Stable stationary solutions with a single internal transition layer when
the bistable nonlinearity is given by the cubic function (\ref{cubic})
with $\alpha = 0.2$. 
The value of the conserved mass $\xi$ is given by $\xi = 0.35$.
Panels (a) and (b) show the graphs of $f(u,v) =0$ for
$\beta = 0.5$ and $\beta = 0.3$, respectively.
Panels (c) and (d) show the profiles of single transition layer solutions for $\beta = 0.5$ and $\beta = 0.3$, respectively,
where the red and blue solid lines indicate the $u$- and 
$v$-components, respectively.
The upper panels of (c) and (d) show that the $v$-components of the stationary solutions appear to be spatially homogeneous.
However, the lower panels, 
enlarged views of the $v$-components show that 
they exhibit a single internal transition layer at $x= x^*$.
The values of $v^*$ defined by $J(v) =0$ are 
given by $v^* = 0.0$ for (c) and $v^* \approx -0.164$ for (d),
respectively, which are indicated by dotted lines.
The layer positions are given by $x^* = 0.650$ for (c)
and $x^* \approx 0.422$ for (d), respectively, which are 
indicated by broken lines. 
These values can be obtained by \eqref{a7} when $\xi = 0.35$.}
\label{example1}
\end{figure}

\begin{example}\label{exa1} \rm

The first example is given by a simple cubic function
\begin{equation}\label{cubic}
f(u, v) = \alpha v - u(u-\beta)(u-1),
\end{equation}
where $\alpha$ and $\beta$ are constants satisfying $\alpha >0$ and $0 < \beta < 1$. 
\end{example}

When $f(u, v) = 0$, we have
\[
v = \frac{1}{\alpha}u(u-1)(u-\beta) =:g(u)
\ \ \
\text{and}
\ \ \ 
\frac{dv}{du} = - \frac{f_u}{\alpha},
\]
where $f_u = -3u^2 + 2(\beta+1)u - \beta$. It is easy to see that the assumptions (A1) and (A3) hold by  $f_v = \alpha$. The ODE $u_t = f(u,v)$ is bistable in $u$ for each 
$v \in (\un{v}, {\ov{v}})$ with $\un{v} = g(h^+(\un{v}))$ and $\ov{v} = g(h^-(\ov{v}))$
as seen in Figure~\ref{example1}(a) and (b). Here, $h^+(\un{v}) = (1 + \beta + \sqrt{\rho(\beta)})/3$
and $h^-(\ov{v}) = (1 + \beta - \sqrt{\rho(\beta)})/3$
are the solutions of $f_u=0$, where $\rho(\beta) = 1 - \beta + \beta^2$.

We immediately find that 
$$
\begin{array}{l}
J'(v) =  \dis\int_{h^-(v)}^{h^+(v)} f_v(u, v)du =  \alpha (h^+(v) - h^-(v)) > 0
\end{array} 
$$
for all $v \in I$. Moreover, we have
$$
J(0) = \int_{h^-(0)}^{h^+(0)} f(u, 0)du  = - \int_0^1 u(u-1)(u-\beta) du = \frac{1}{12} - \frac{\beta}{6} 
$$
by $h^-(0) = 0$ and $h^+(0) = 1$. 

When $\beta = 1/2$, it is easy to see that the assumption (A2) 
holds because $J(v^*)= 0$ and $J'(v^*) > 0 $ when $v^* =0$.
In addition, by a formal matched asymptotic expansion, 
we can give an approximate expression of a 
stationary solution with
a single internal transition layer
$$
u^\ep(x) = \dfrac{1}{2} \left\{ 1 + \tanh \left( \dfrac{x - x^*}{2 \sqrt{2} \ep} \right) \right\} + O(\ep)
$$
and 
$$ 
v^\ep(x) = v^* - \dfrac{\ep^2}{2D}\tanh \left( \dfrac{x - x^*}{2 \sqrt{2} \ep} \right) + O(\ep^3), 
$$
where the layer position $x^* = 1 - \xi$ is given by \eqref{a7}.
The leading terms of $u^\ep(x)$ and $v^\ep(x)$ characterize 
the profile of the stationary solution in 
the upper panel in Figure~\ref{example1}(c). 
The $O(\ep^2)$ term of $v^\ep(x)$
characterizes the transition layer of the $v$-component of
the stationary solution in the lower panel in Figure~\ref{example1}(c). This $O(\ep^2)$ term is necessary to 
precisely obtain the profile of the stationary solution, whereas
it is not necessary to determine its stability as 
seen in Section 3.

When $\beta \neq 1/2$, we can factorize 
$f(u, \un{v})$ and $f(u, \ov{v})$ into 
$
f(u, \un{v})= (u - h^+(\un{v}))^2 (h^-(\un{v}) - u)
$
and 
$
f(u, \ov{v})= (u - h^-(\ov{v}))^2 (h^+(\ov{v}) - u),
$
respectively, where
$$ 
\displaystyle{ h^-(\un{v}) = \frac{1 + \beta - 2 \sqrt{\rho(\beta)} }{3} }
\ \ \text{and} \ \ 
\displaystyle{ h^+(\ov{v}) =  \frac{1 + \beta + 2 \sqrt{\rho(\beta)} }{3} }. 
$$
Then, we have
$$
J (\un{v})  =  \int_{h^-(\un{v})}^{h^+(\un{v})} f(u, \un{v})du  = - \frac{\rho(\beta)^2}{12} < 0
$$
and 
$$ 
J (\ov{v})  =  \int_{h^-(\ov{v})}^{h^+(\ov{v})} f(u, \ov{v})du =  \frac{\rho(\beta)^2}{12} > 0. 
$$
Therefore, we find that $J(v^*) = 0$ and $J'(v^*) > 0$ hold for some $v^* \in I$, which implies that the assumption (A2) holds for $\beta \neq 1/2$.
Figure \ref{example1}(d) shows a stable stationary solution with a single internal 
transition layer, in which the value of $v^*$ appears in $(\un{v}, 0)$ due to 
$\beta \neq 1/2$.

\begin{figure}[ht!]
 \centering
\includegraphics[width=12cm, bb=0 0 890 571]{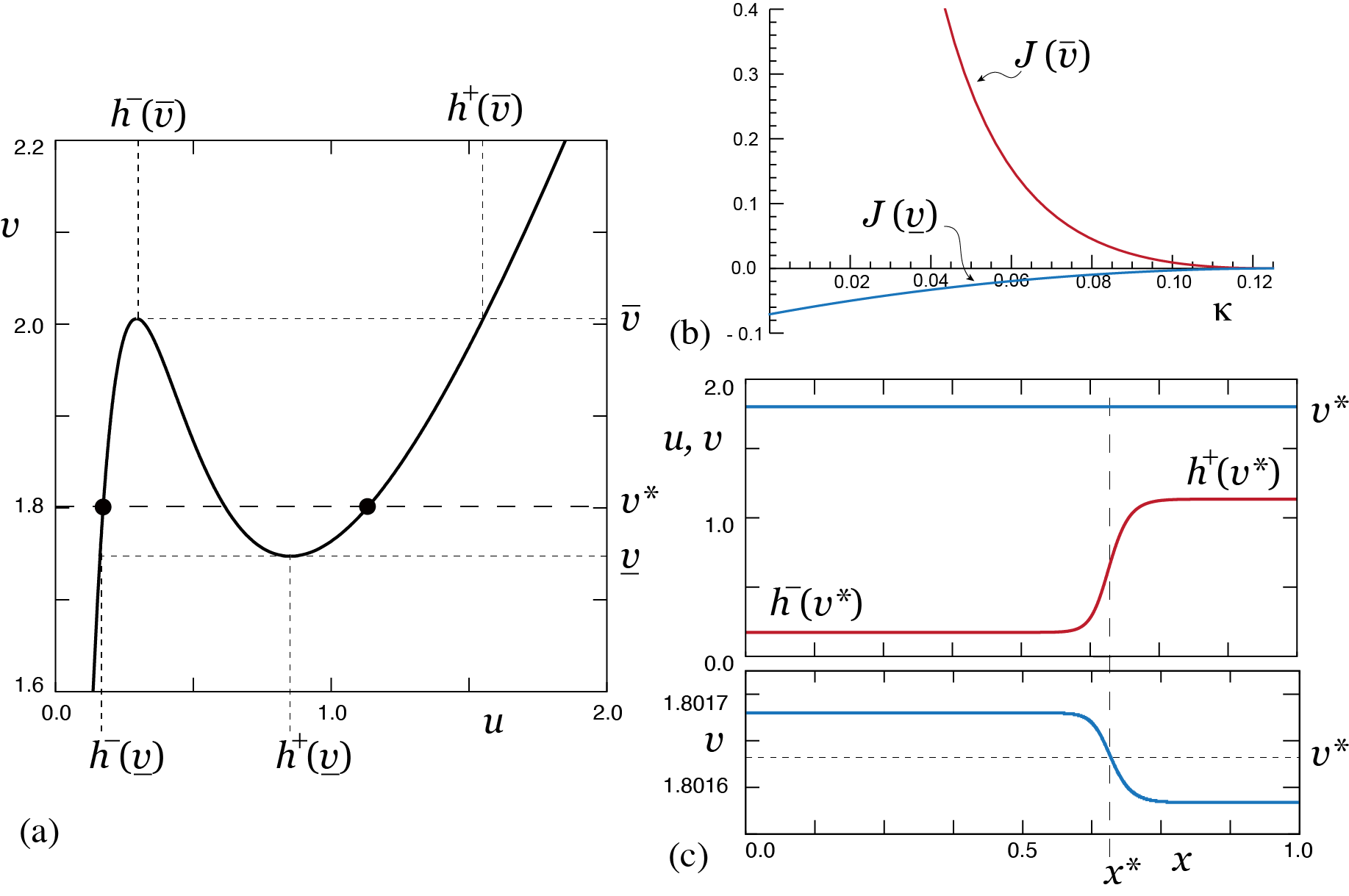}
  \caption {Stable stationary solutions with a single internal transition layer when
the bistable nonlinearity is given by the Hill type function \eqref{mori}. 
The value of the conserved mass $\xi$ is given by $\xi = 2.3$.
(a) The graph of $f(u,v) =0$ for $\kappa = 0.067$. 
(b) The graphs of $J (\un{v})$ and $J (\ov{v})$, which are continuous functions in $\kappa$.  
(c)  The profile of a single transition layer solution; the red and blue solid lines indicate the $u$- and $v$-components, respectively. 
The lower panel of (c),
an enlarged view of the $v$-component shows that 
it exhibits a single internal transition layer at $x= x^*$.
The value of $v^*$ defined by $J(v) = 0$ is given by $v^* \approx
1.802$.
The layer position is given by $x^* \approx 0.660$, which  
can be obtained by \eqref{a7} when $\xi = 2.3$.
}
\label{example2}
\end{figure}
\begin{example}\label{exa2}  \rm
The second example is given by a Hill type function
\begin{equation}\label{mori}
f(u, v) = 
\displaystyle{
\left( \kappa + \frac{u^2}{1+u^2} \right) v - u,
}
\end{equation}
where $\kappa$ is a positive constant satisfying $0 < \kappa < 1/8$.
\end{example}

This model proposed by \cite{MJE2} has a bistable nonlinearity. In fact, 
we can numerically obtain the curve of $f(u,v) = 0$ for $\kappa = 0.067$, 
as shown in Figure~\ref{example2}(a). 
The ODE $u_t = f(u, v)$  is bistable in $u$ for $v \in I = (\un{v}, \ov{v})$ with  
$$
\displaystyle{ \un{v} = \frac{1 + \kappa - \phi}{2 \kappa (1+\kappa) }  \sqrt{ \frac{\phi}{1+ \kappa} } }
\ \ \text{and} \ \  
\displaystyle{ \ov{v} = \frac{1 + \kappa - \psi}{2 \kappa (1+\kappa) }
 \sqrt{ \frac{\psi}{1+ \kappa} } },
$$
where 
$$
\displaystyle{  \phi(\kappa) = \frac{1 - 2 \kappa + \sqrt{1-8 \kappa} }{2}  }
\ \ \text{and} \ \  
\displaystyle{  \psi(\kappa) = \frac{1 - 2 \kappa - \sqrt{1-8 \kappa} }{2}  } .
$$
Moreover, with the aid of the software MATHEMATICA, we can obtain
$$
\displaystyle{ h^+(\un{v}) =  \sqrt{ \frac{\phi}{1+ \kappa} }}
\ \ \text{and} \ \ 
\displaystyle{ h^-(\un{v}) = \frac{1 + \kappa - \phi }{2 \sqrt{(1 + \kappa) \phi}} }
$$
and
$$
\displaystyle{ h^+(\ov{v}) = \frac{1 + \kappa - \psi }{2 \sqrt{(1 + \kappa) \psi}} }
\ \ \text{and} \ \ 
\displaystyle{ h^-(\ov{v}) =  \sqrt{ \frac{\psi}{1+ \kappa} }}.
$$
Furthermore, we have
\begin{eqnarray*}
\displaystyle{ \frac{d h^-(\un{v})}{d \kappa} } & = & 
\displaystyle{\frac{(1+\kappa + \phi)(\phi - (1+\kappa) \phi'(\kappa))}{4 (\sqrt{(1+\kappa) \phi})^3}} > 0
\end{eqnarray*}
for $0 < \kappa < 1/ 8$, which implies that $h^-(\un{v})$ is 
a monotone increasing function in $\kappa \in (0, 1/8)$.
Since $\lim_{\kappa \to 0}  h^-(\un{v}) = 0$ 
and $\lim_{\kappa \to 1/8}  h^-(\un{v}) = 1/ \sqrt{3}$, 
we have $0 < h^-(\un{v}) < 1/\sqrt{3}$. Similarly, we can show that $h^+(\ov{v})$ is a monotone decreasing function in $\kappa \in (0, 1/8)$ with $1/\sqrt{3} < h^+(\ov{v}) < +\infty $. 
Therefore, noting
$$
\displaystyle{f_u(u, v) = \frac{2 u v}{(1+u^2)^2}-1}
\ \ \text{and} \ \ 
\displaystyle{f_v(u, v) = \kappa + \frac{u^2}{1+u^2}},
$$
we see that the assumptions (A1) and (A3) hold.

Next, we check the assumption (A2). A direct calculation shows
$$
\begin{array}{rcl}
J'(v) & = &
\displaystyle{ \int_{h^-(v)}^{h^+(v)} f_v(u, v)du  } 
\\[2ex]
& = &
\dis{(h^+(v) - h^-(v)) 
\left( \kappa + 1 - \frac{\tan^{-1} (h^+(v)) - \tan^{-1} (h^-(v)) } {h^+(v) - h^-(v)} \right) }  
\\[3ex]
 & > &  \kappa (h^+(v) - h^-(v)),
\end{array}
$$
which implies that $J'(v) > 0$ holds for all $v \in I$. Moreover, we have 
$$
J (\un{v})  =  
\frac{(\kappa - 2 \phi)(1 + \kappa - \phi)^2 + 4 \phi^2 (1 - \phi) }{8 \phi \kappa (1 + \kappa) }
- \frac{1 + \kappa - \phi}{2 \kappa (1 + \kappa)} \sqrt{ \frac{\phi}{1 + \kappa}  }
\cdot K_1
$$
and 
$$ 
J (\ov{v})  =  
- \frac{(\kappa - 2 \psi)(1 + \kappa - \psi)^2 + 4 \psi^2 (1 - \psi) }{8 \psi \kappa (1 + \kappa) } 
  +   \frac{1 + \kappa - \psi}{2 \kappa (1 + \kappa)} \sqrt{ \frac{\psi}{1 + \kappa}  }
\cdot K_2,
$$
where
$$
K_1 = 
\tan^{-1} \left(  \sqrt{ \frac{\phi}{1 + \kappa} } \right) - 
\tan^{-1} \left(  \frac{1 + \kappa - \phi}{2 (1 + \kappa) }  \sqrt{ \frac{1 + \kappa}{\phi} } \right)   
$$
and
$$ 
K_2 = 
\tan^{-1} \left(  \sqrt{ \frac{\psi}{1 + \kappa} }  \right) -
\tan^{-1} \left(  \frac{1 + \kappa - \psi}{2 (1 + \kappa) }  \sqrt{ \frac{1 + \kappa}{\psi} } \right). 
$$
As shown in Figure~\ref{example2}(b), we see that 
$J (\un{v})$ and $J (\ov{v})$ are continuous functions 
in $\kappa$, and that
$J (\un{v}) < 0$ and $J (\ov{v}) > 0$ for any $\kappa \in (0, 1/8)$.
In addition, we have
$\lim_{\kappa \to 1/8} J (\un{v}) =  \lim_{\kappa \to 1/8} J (\ov{v}) = 0$, 
$\lim_{\kappa \to 0} J (\un{v}) = (3 - \pi)/2 < 0 $ and $\lim_{\kappa \to 0} J (\ov{v}) = + \infty$. 
Therefore, we find that $J(v^*) = 0$ and $J'(v^*) > 0$ hold for some $v^* \in I$, which implies
that the assumption (A2) holds.

Figure~\ref{example2}(c) shows a stable stationary solution with a single internal 
transition layer, in which the layer position $x^*$ given by \eqref{a7} is numerically computed.
 \begin{figure}[ht!]
 \centering
\includegraphics[width=12cm, bb=0 0 900 726]{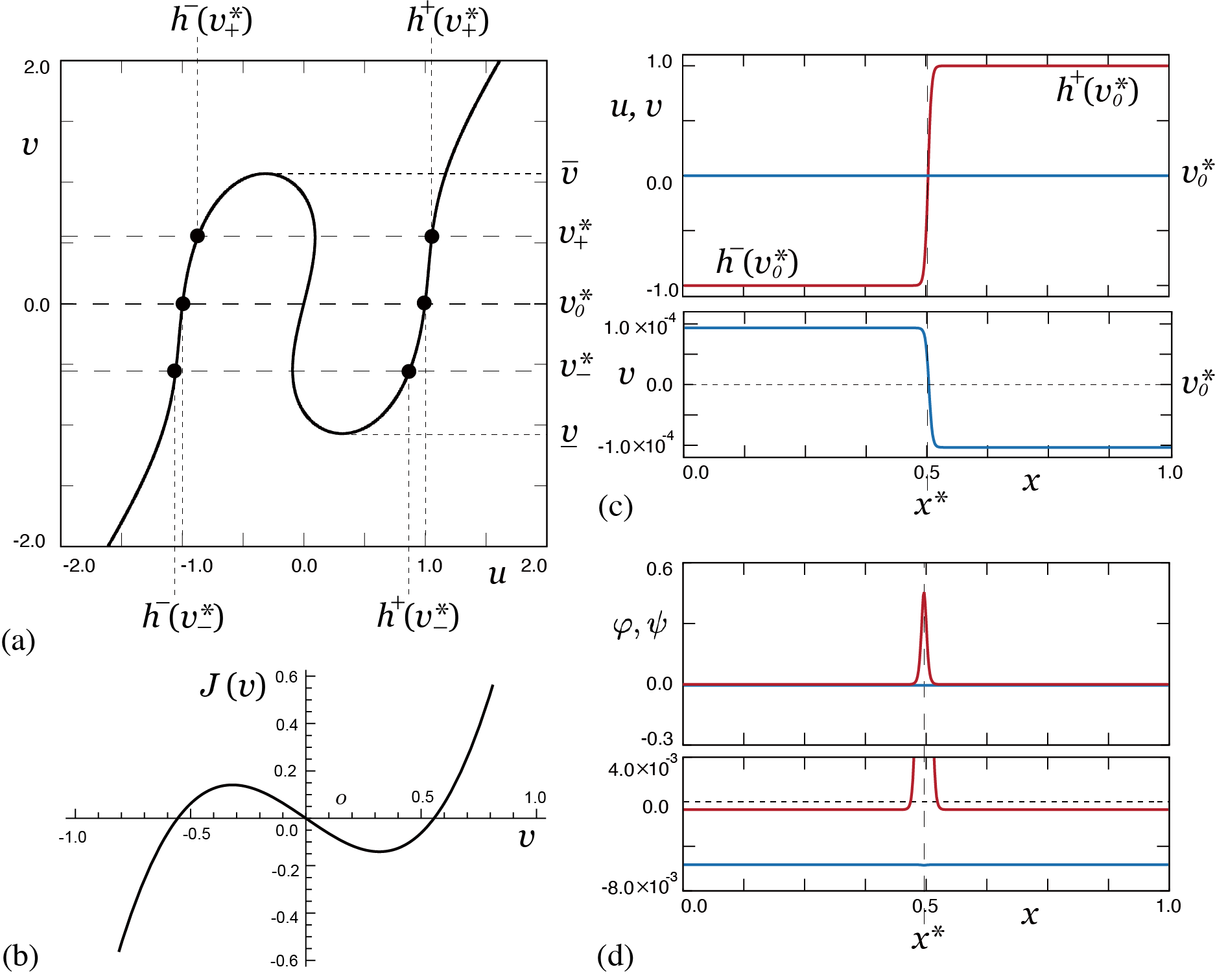}
  \caption {Unstable stationary solutions with a single internal transition layer when the bistable nonlinearity is given by
the artificial function \eqref{multivalued}. 
The value of the conserved mass $\xi$ is given by $\xi = 0.0$.
Panels (a) and (b) show a curve defined by $f(u,v) =0$ and the graph of $J(v)$ defined by \eqref{a4}, respectively. In (a), $\ov{v} = -\un{v} \approx 1.072$. The solid points indicate the bistable equilibria of the ODE $u_t = f(u,v)$ in 
$u$ for $v= v^*_0$ and $v = v^*_\pm$, which are three zeros of $J(v)$ in (b). The values $v^*_0$,  $v^*_+$ and $v^*_-$ are given by
$v^*_0 = 0.0$,  $v^*_+ \approx 0.55589$ and $v^*_- \approx -0.55589$, respectively.
Panel (c) shows the profile of an unstable stationary solution with a single internal transition layer, where
the red and blue solid lines indicate the $u$- and $v$-components, respectively.  
The lower panel of (c),
an enlarged view of the $v$-component shows that
it exhibits a single internal transition layer at $x = x^*$.
The layer position is given by $x^* = 0.5$, which  
can be obtained by \eqref{a7} when $\xi = 0.0$.
Panel (d) shows the profile of the eigenfunction $\Phi^{\ep} = (\varphi^\ep, \psi^\ep)$ associated with the critical eigenvalue
$\lambda^{\ep} \approx 6.2 \times 10^{-3}$ characterized by \eqref{cri}, where
the red and blue solid lines indicate $\varphi^\ep$ 
and $\psi^\ep$, respectively. 
The lower panel of (d), an
enlarged view of the profile of $\Phi^{\ep}$ is useful to understand
that $\Phi^{\ep}$ satisfies the constrained condition in \eqref{c2}.
}
 \label{example3}
\end{figure}

\subsection{Unstable case}

In contrast to the stable case, it is not easy to
give a mass-conserving reaction-diffusion model which has
unstable stationary solutions with a single internal transition layer.
To the best of our knowledge, such a model 
has not yet been proposed. In fact, 
a curve defined by $f(u,v)=0$ cannot
be expressed by a function of $u$ under Assumption~\ref{ass1}
when $J'(v^*) < 0$.
The same problem has also been well-recognized in the case of 
reaction-diffusion systems of FitzHugh-Nagumo
type. According to \cite{NF}, such a reaction-diffusion 
system can have unstable transition layer solutions if 
the null cline of the reaction term with cubic nonlinearity is given by a curve 
as shown in \cite[Figure 11]{NF}.
However, \cite{NF} presented no 
models with such nonlinear reaction terms; they have not yet 
been proposed as well as the case of mass-conserving 
reaction-diffusion systems with bistable nonlinearity.

\begin{example}\label{exa3}  \rm

The third example is given by an artificial function
\begin{equation}\label{multivalued}
f(u, v) = 4 u - v - 4 u^3 + 2 u^2 v - \frac{5}{4} u v^2 - \frac{3}{4} v^3,
\end{equation}
which satisfies $f(-u, -v) = - f(u, v)$ and $f(\pm 1, 0) = 0$. 
\end{example}

Figure~\ref{example3}(a) shows a curve defined by $f(u,v) = 0$, 
which is similar to that in \cite[Figure 11]{NF}. 
We can verify that the assumptions (A1) and (A3) hold because
the ODE $u_t = f(u,v)$ is bistable in $u$ for each $v \in (\un{v}, \ov{v})$, 
and the values of $\ov{v}$ and $\un{v}$ are determined by $f(u,v)=0$ and $f_u(u,v) = 0$. 
Moreover, we can obtain the graph of $J(v)$ shown
in Figure~\ref{example3}(b)
with the aid of the software MATHEMATICA.
The values $v^*_0$,  $v^*_+$ and $v^*_-$ defined by
$J(v) =0$ satisfy $J'(v^*_0) < 0$ and $J'(v^*_\pm) > 0$,
which implies that the assumption (A2) holds.
It is easy to see that $v^*_0 = 0$, $h^{+}(v^*_0) = 1$ and $h^{-}(v^*_0) = -1$.
The values of $v^*_\pm$, $h^{+}(v^*_\pm)$ and $h^-(v^*_\pm)$
are numerically obtained as
$v^*_+ = - v^*_- \approx 0.55589$, $h^{+}(v^*_+) = - h^{-}(v^*_-) \approx 1.05659$ and $h^{+}(v^*_-) = - h^{-}(v^*_+) \approx 0.87129$.

The characterization of the critical eigenvalue given by 
\eqref{cri} suggests that the sign of $J'(v)$ at the zero of $J(v)$ can determine the
stability/instability of a single internal transition layer solution. 
We can obtain an unstable stationary solution with a single internal transition layer
as shown in Figure~\ref{example3}(c)
by applying the Newton method to \eqref{b1};
the initial value of the iteration is given by a pair of functions obtained from
$(U^*(x), V^*(x))$ by replacing $v^*$ with $v^*_0$, where $(U^*(x), V^*(x))$
is given by \eqref{a8} and \eqref{a9}.
The layer position is given by $x^* = 0.5$, which can be
obtained by \eqref{a7} when $\xi = 0.0$.
Moreover, by numerically solving the linearized eigenvalue problem \eqref{c1} around the stationary solution, 
we find that the critical eigenvalue characterized by \eqref{cri} is given by $\lambda^\ep \approx 6.2 \times 10^{-3} > 0 $, and the profile of 
its associated eigenfunction $\Phi^\ep = (\varphi^\ep, \psi^\ep) \in X$
is shown in Figure~\ref{example3}(d), where $X$ is given by \eqref{c2}.

\begin{figure}[ht!]
 \centering
\includegraphics[width=12cm, bb=0 0 905 755]{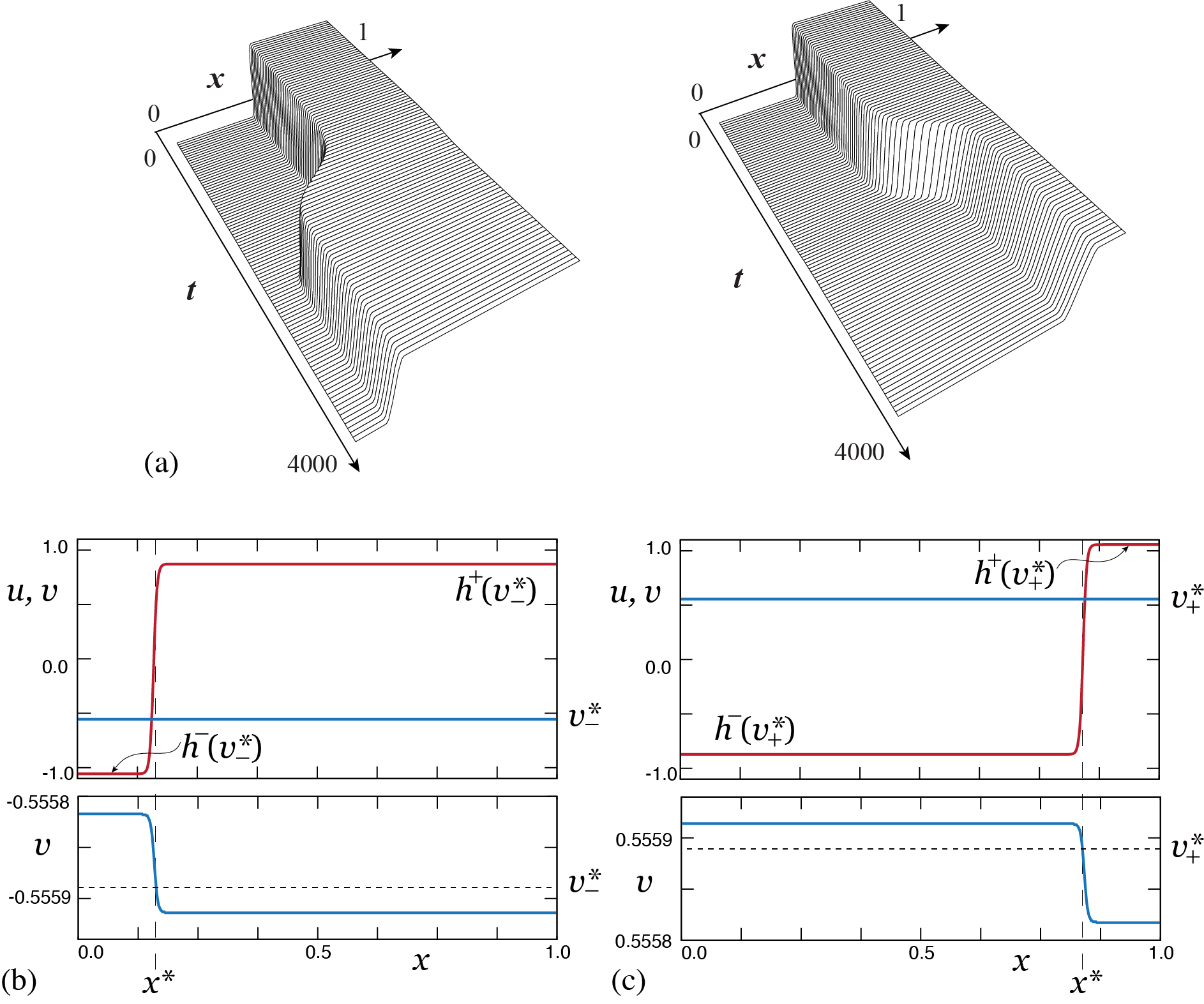}
  \caption {Panel (a) shows that solutions starting from 
the unstable stationary solution in Figure~\ref{example3}(c) with small perturbations
converge to one or the other of
stable transition layer solutions corresponding to 
$v^*_-$ and $v^*_+$ defined by $J(v) = 0$.
Panels (b) and (c) show the profiles of these stable solutions when the bistable nonlinearity is given by \eqref{multivalued}, 
where
the red and blue solid lines indicate the $u$- and $v$-components, respectively.  
The lower panels of (b) and (c) are
enlarged views of the profiles of the $v$-components. The value of the conserved mass $\xi$ is given by $\xi = 0.0$. 
The positions of the transition layers in (b) and (c) are $x^* \approx 0.164$ and $x^* \approx 0.836$, respectively, which are indicated by broken lines. }
 \label{example4}
\end{figure}

When we add small perturbations
to the unstable stationary solution in the direction of 
the associated eigenfunction $\Phi^\ep$, 
numerical solutions converge to one or the other of
stable transition layer solutions shown in Figure~\ref{example4}(b)
and (c). These stable solutions are corresponding to 
$v^*_-$ and $v^*_+$ defined by $J(v) = 0$.
We note that their stability is guaranteed by Theorem~\ref{th4} because $J'(v^*_\pm) > 0$.
The positions of the transition layers in 
Figure~\ref{example4}(b) and (c) are given by $x^* \approx 0.164$ and $x^* \approx 0.836$, respectively. These values can be obtained by 
\eqref{a7} when $\xi = 0.0$.
 
\begin{remark} \rm
In general, from the topology of a curve defined by $f(u,v)=0$ under Assumtion~\ref{ass1},
if there exists a
zero of $J(v)$ with $J'(v) < 0$, then there exists two zeros of $J(v)$ with $J'(v) > 0$.
This implies that if there exists an unstable transition layer solution, then
there exist two stable transition layer solutions. 
We note that this unstable transition layer solution is a separatrix of the dynamics of 
mass-conserving reaction-diffusion systems with bistable nonlinearity.
\end{remark}

\begin{remark} \rm
Figure~\ref{example3}(d) suggests that 
$|| \psi^\ep_x ||_{L^2} = O(\ep^p)$ for some $p >0$.
This information may be useful to prove the existence of
the critical eigenvalue, which leads to the instability of 
single transition layer solutions when $J'(v^*) < 0$.
\end{remark}

\section{Concluding remarks} \label{conclude}

In this paper, we show that mass-conserving reaction-diffusion systems with 
bistable nonlinearity can have stable stationary solutions 
with a single internal transition layer under general assumptions.
Our approach is based on the singular perturbation method by \cite{HS, I, MTH, NF}.
In spite of the complication that we cannot apply 
the Lax-Milgram theorem to the singular limit eigenvalue problem (SLEP) concerning
the stability of the stationary solutions, 
we can give the precise characterization of a critical eigenvalue 
which essentially determines the stability 
due to a natural constrained condition derived from the conservation law.
Consequently, we can 
provide a rigorous proof of the stability of the stationary solutions under general 
assumptions. 
Our result provides a basic information for studying the bifurcation structure of
mass-conserving reaction-diffusion systems with bistable nonlinearity
\cite{MKNTY, MJE2}.
Moreover, we emphasize that our approach can apply to singular perturbation problems in higher spatial dimensions as seen in 
reaction-diffusion systems with bistable nonlinearity of FitzHugh-Nagumo 
type \cite{SS1, SS2}.

Our characterization of the critical eigenvalue suggests that the stationary solutions 
with a single internal transition layer are unstable if $J'(v^*) < 0$.
That is, the sign of $J'(v^*)$ determines
the stability/instability of the stationary solutions; 
it is a conjecture for future studies.
Since the spatial dimension of our problem is one, 
we can apply a different method by \cite{I}
for studying the stability of the stationary solutions, which 
is substantially based on the ideas of the Evans function theory. In this case, 
we also encounter a problem related to that in Remark~\ref{rem3y}, which would require a delicate mathematical analysis. 
At present, we can prove that this conjecture is true if we replace the assumption (A3) 
in Assumption~\ref{ass1} by a somewhat stronger one.

We can construct stationary solutions with multiple internal 
transition layers following a standard procedure in \cite{N}.
For example, we can construct a stationary solution with
double layers $(u^\ep(x), v^\ep(x))$ as 
follows: Let $(\tilde{u}^\ep(x), \tilde{v}^\ep(x))$ be
a stationary solution with a single layer on $[0, 1/2]$, which is 
constructed by our method. Then,
the values of $(u^\ep(x), v^\ep(x))$ for $x \in [0, 1/2]$ 
and for $x \in [1/2, 1]$ are defined 
by $(\tilde{u}^\ep(x), \tilde{v}^\ep(x))$
for $x \in [0, 1/2]$ and $(\tilde{u}^\ep(1-x), \tilde{v}^\ep(1-x))$
for $x \in [1/2, 1]$, respectively.
However, the stability analysis for multi-layered stationary solutions seems quite 
delicate. In fact, mentioned in \cite[Section 3.3]{MJE2}, 
multi-layered solutions may move very slowly at the speed of $O(\ep^{-C/\ep})$.
This suggests that 
the linearized operator around 
a multi-layered stationary solution can have
$O(e^{-C/\ep})$-eigenvalues, and it is not easy to determine
the sign of these extremely small eigenvalues.
 
We can consider simple models obtained by adding 
perturbations to mass-conserving reaction-diffusion systems with 
bistable nonlinearity. In fact, \cite{VC} studied such a model
obtained by adding source and loss terms to a mass-conserving
reaction-diffusion system
for studying the wave-pinning phenomenon in cell division and 
differentiation \cite{MJE2}.
These terms break the mass conservation law, and 
give rise to several equilibrium solutions such as
homogeneous equilibria, periodic patterns, localized patterns,
and isolated spikes. It is an interesting topic to 
investigate the existence and stability of
transition layer solutions in such perturbed systems
without mass conservation. This topic can provide  
a viewpoint for understanding a transition from 
localized patterns to layered patterns in reaction-diffusion systems, 
which was studied by (formal) matched asymptotic expansions and numerical 
simulations \cite{VC}.

We consider mass-conserving reaction-diffusion systems with bistable nonlinearity,
in which the ODEs $u_t = f(u, v)$ 
obtained by dropping the diffusion terms are bistable in $u$ for each fixed $v$. 
They are related to the wave-pinning phenomenon in cell division and 
differentiation \cite{MJE1,MJE2}. 
On the other hand, 
we can consider other mass-conserving reaction-diffusion systems with bistable nonlinearity,
in which the ODEs $v_t = - f(u, v)$ are bistable in $v$ for each fixed $u$. 
It was shown in \cite{KIE, OkS} that these systems are related to 
cell polarity oscillations, which cause the reversal of cell polarity. 
We expect that mathematical analyses for mass-conserving reaction-diffusion systems with 
bistable nonlinearity can provide a significant viewpoint for 
understanding cell polarity formation, which plays a key role in cell division and differentiation.

\vspace{1cm}

{\bf Acknowledgments.} 
The authors would like to express their appreciation to the referees for their useful suggestions and comments, which have improved the original manuscript.
The first and second authors were supported 
in part by the JSPS KAKENHI Grant Numbers JP24K06845 and 
JP23K03209, respectively. The third author was also supported 
in part by the JSPS KAKENHI Grant Numbers JP19K03618 and JP24K06864.

\end{document}